%% file: gaspin.tex
\documentclass[final]{siamltex}

\usepackage{graphicx,multicol}        % used for the two-column index
\usepackage[bottom]{footmisc}% places footnotes at page bottom
\usepackage{amsmath,amssymb}

\usepackage{algorithm}

\usepackage{algpseudocode}

\newcommand{\eref}[1]{Equation~\eqref{#1}}
\newcommand{\D}{{\cal S}}

\newcommand{\ATR}{A${}_{\text{tr}}$}
\newcommand{\pred}{{\mathop{\mathrm{pred}}\nolimits}}
\newcommand{\ared}{{\mathop{\mathrm{ared}}\nolimits}}

\algnewcommand\algorithmicinput{\textbf{Input:}}
\algnewcommand\Input{\item[\algorithmicinput]}

\algnewcommand\algorithmicoutput{\textbf{Output:}}
\algnewcommand\Output{\item[\algorithmicoutput]}

\newtheorem{thm}{Theorem}[section]

\newtheorem{lem}[thm]{Lemma}
\newtheorem{remark}[thm]{Remark}

\usepackage{todonotes}

\usepackage{bibentry}

\input{commands.tex}

\title{On the Globalization of ASPIN employing Trust-Region Control Strategies -- 
	Convergence Analysis and Numerical Examples}

\author{Christian Gross\thanks{Institute of Computational Science and 
Euler Institute,
Universit\`a della Svizzera italiana (USI), 
via G. Buffi 13, 6900 Lugano, Switzerland (\tt{rolf.krause@usi.ch})},
\and Rolf Krause${}^*$ }

\begin{document}

\nobibliography*

\maketitle                   % Produces the title.

\begin{abstract}
  The parallel solution of large scale non-linear programming problems, which arise for example from the discretization of non-linear partial differential equations, is a highly demanding
task. Here, a novel solution strategy is presented, which is inherently parallel and globally convergent. Each global non-linear iteration step consists of
asynchronous solutions of local non-linear programming problems followed by a global recombination step. The recombination step, which is the solution of a
quadratic programming problem, is designed in a way such that it ensures global convergence. 
As it turns out, the new strategy can be considered as a globalized additively preconditioned inexact Newton (ASPIN) method~[\nobibentry{CaiKeyes00}]. 
However, in our approach the influence of ASPIN's non-linear preconditioner on the gradient  is controlled in order to ensure a
sufficient decrease condition. Two different control strategies  are described and analyzed.  Convergence to first-order critical points of our non-linear
solution strategy is shown under standard trust-region assumptions.  The strategy is investigated along difficult minimization problems arising from non-linear
elasticity in $3D$ solved on a massively parallel computer with several thousand cores.
\end{abstract}

\section{Introduction}

The massive use of parallel computers with hundreds and thousands 
of processors enforces algorithms to be designed especially for parallel computing. 
This means that large scale problems must be divided into 
smaller subproblems in order to be solvable on modern super computers. 
Here, we consider a novel decomposition approach for solving the following
smooth non-linear programming problem
\begin{equation}
\label{eqn:nlMinProb}
	u\in \mathbb R^n: J(u) = \min!,
\end{equation}
where $J:\mathbb R^n\to \mathbb R$ is a continuously differentiable objective
function. 

Since the objective function in general might be non-convex, one has to employ a \emph{globalization strategy} such as trust-region methods 
(for a broad introduction see \cite{ConnGouldToint00}) or linesearch strategies (for an overview see \cite{NocedalWright2006}) 
in order to ensure the convergence to a local minimizer of \eref{eqn:nlMinProb}. 
An important feature of these strategies is that the way corrections are computed is arbitrary, as long as what is known as a sufficient decrease condition is satisfied.  
A straightforward approach for parallelizing these strategies is to parallelize the computation of the corrections and the assembling process.

As an alternative -- and possibly in a more parallel spirit --, the global problem \eqref{eqn:nlMinProb} can be solved by splitting it into local
non-linear subproblems, which are then solved asynchronously and in parallel. In this field, several different classes of inherently parallel
globalization strategies were developed, such as the parallel variable and gradient distribution~\cite{FerrisMangasarian1994,Mangasarian95} and the additively
preconditioned trust-region (APTS) and linesearch methods~\cite{CGross_2009}. 
For a recent parallel method for constrained minimization problems see for example \cite{RKrause_ARigazzi_and_JSteiner_2015}.

In~\cite{CaiKeyes00}, a non-linear left preconditioned solution strategy called additively preconditioned inexact Newton (ASPIN) method 
was presented. This meth\-od is designed for the solution of non-linear equations of the kind
\begin{equation}
\label{eqn:nlEq}
	u\in\mathbb R^n: F(u) = 0,
\end{equation}
where $F:\mathbb R^n \to \mathbb R^n$. Obviously, the first-order condition of \eref{eqn:nlMinProb} 
\begin{equation}
\label{eqn:nlFirstOrder}
	u\in\mathbb R^n: \nabla J(u) = 0
\end{equation}
is of the kind introduced in \eref{eqn:nlEq}. 
Although ASPIN is designed as an inherently parallel non-linear solution strategy, it does not claim to be a globalization strategy. 
It is particularly difficult to prove global convergence of this strategy, 
since the computation of the global correction is intertwined with the solution of the local non-linear programming problems,
as we will see in Section~\ref{sec:aspin}. 
Thus, only local convergence results are available, 
i.e., it was shown that ASPIN solves \eref{eqn:nlEq} if the initial iterate is sufficiently close to a local solution, see e.g.~\cite{CaiKeyes00,ArnalEtAl2008}.  
Despite the lack of global convergence results, 
ASPIN has been successfully employed as a solution strategy for non-linear problems arising for instance in the field of fluid dynamics~\cite{CaiKeyes2002a,CaiKeyes2002}.

Here, we propose a new family of non-linear and parallel solution strategies, which are globally convergent and which -- near a minimizer -- reduce to the ASPIN method. In this sense, the members of this family can be viewed as globalized ASPIN methods. As we show in our analysis, ASPIN's Newton corrections can be considered to be the solutions of the first-order conditions of perturbed or \emph{non-linearly preconditioned} quadratic programming problems. The arising sequence of quadratic programming problems then can be embedded into a trust-region framework, thus allowing for global convergence control. 
However, within the trust-region framework the influence of the non-linear preconditioning step has to be  taken into account.  Therefore, we accompany the trust-region control by  a criterion for controlling the influence of the non-linear preconditioning on the quadratic programming problem.

This criterion is inspired from the analysis of trust-region methods with perturbed models, where inexactly computed gradients occur.  
In the literature, two different approaches  can be found for handling the perturbation of the gradient. 
Usually, the perturbation is assumed to be bounded either by the norm of the exact gradient (cf., \cite{Carter1993}) or by the trust-region radius itself (cf., \cite{Toint1988,ConnEtAl1993}). 
Here, we exploit the second approach and introduce a second, local trust-region radius. This second trust-region radius is then employed for controlling the parallel local solution phase. 
In particular, we present two different approaches, which lead to preconditioned gradients satisfying the control criterion: 
a modified, local trust-region method similar to the local solution phase of the APTS method~\cite{CGross_2009}
and a damping approach.

An important result of our analysis is that standard trust-region assumptions, which in particular do not include convexity assumptions, allow to prove convergence of the new ``G-ASPIN'' method, cf.~Section~\ref{sec:firstOrderConvergence}. On the other hand, the numerical studies in Section~\ref{sec:numerics} show that our globalized ASPIN method efficiently resolves local non-linearities and yields fast convergence for large-scale non-linear programming problems.

\section{The Domain Decomposition Framework}

Minimization problems of the kind \eqref{eqn:nlMinProb} often arise from non-linear programming problems stated in a finite dimensional space $\mathcal X$,  
for instance a finite element space. 
Then a coordinate isomorphism $X:\mathbb R^n\to \mathcal X$  exists, 
which maps coefficient vectors from $\mathbb R^n$ to elements in $\mathcal X$. 
We assume that $\mathcal X$ is decomposed into $N$ subspaces $\mathcal X^1,\ldots,\mathcal X^N$ and
obtain that also $\mathbb R^n$ is decomposed into $N$ subspaces,
which will be denoted by $\Sk \subseteq \mathbb R^n$ for $k=1,\ldots,N$. 

In this case, we can define the \emph{prolongation operator} $\Ik:\Sk\to \mathbb R^n$ as
\begin{equation}
\label{eqn:definitionIk}
	X \Ik \localu = \Xk \localu \qquad \text{ for all } \localu\in \Sk
\end{equation}
where $\Xk:\Sk \to \mathcal{X}^k$ is the subset coordinate isomorphism.
As common for example also in multigrid literature, 
the \emph{restriction operator} $\Rk$ is given by $\Rk = (\Ik)^T$.

The restriction operator, as a transfer operator, aims at transferring dual quantities, i.e., the residual, from $\mathbb R^n$ to $\Sk$. 
Furthermore, in non-linear domain decomposition methods, it is necessary to transfer also iterates from $\mathbb R^n$ to $\Sk$. 
To this end, we follow \cite{CGross_RKrause_2008} and introduce a 
\emph{projection operator} $\Pk:\mathbb R^n\to\Sk$ for the projection of primal variables. 
For $u\in \mathbb{R}^n$, the projection~$\Pk u \in \Sk$ is defined as
\begin{equation}
	\label{eqn:defPk}
	\norm{X \left (\Ik \Pk u - u 	\right)}_{\mathcal X} \le 
	\norm{X \left (\Ik v^k - u		\right)}_{\mathcal X} \qquad \text{for all } v^k\in\Sk.
\end{equation}
For further information, in particular, how to practically compute $\Pk$, we refer to~\cite{CGross_2009}. 

Finally, we remark that matrix and vector norms are meant to be $\|\cdot\| = \|\cdot\|_2$.

\section{The ASPIN Method}
\label{sec:aspin}
In~\cite{CaiKeyes00}, X.-C. Cai and D.E. Keyes introduced the concept of 
preconditioned inexact Newton (PIN) methods and in particular ASPIN. 
The ASPIN method is an iterative method for the solution of problems of the kind~\eqref{eqn:nlFirstOrder}. 
For a given iterate $\globiti \in \mathbb R^n$, the new iterate $u_{i+1}$ is defined as
\begin{equation}
	\label{eqn:iterativeStep}
	u_{i+1} = \globiti+\alpha_i \globstepi,
\end{equation}
where $\alpha_i\in(0,1]$ is a damping parameter and $\globstepi\in \mathbb R^n$ is the correction in step~$i$.
The basic principle of ASPIN is to reformulate the computation of the original Newton step~$\globstepi$
\begin{equation}
	\label{eqn:newtonStep}
	\nabla^2 J(\globiti)\globstepi = -g_i
\end{equation}
to the following step
\begin{equation}
\label{eqn:aspinNewtonStep}
	C^{-1}_i \nabla^2 J(\globiti)\globstepi = -g_i^{\text{ASPIN}},
\end{equation}
where $g_i=\nabla J(\globiti)$ is the non-linear residual.
The other two unknown quantities in  \eref{eqn:aspinNewtonStep}, that is the vector $g_i^{\text{ASPIN}}$ and the operator $C^{-1}_i$,
will be explained in detail in the following paragraphs.

We start with the definition of the right-hand side of \eref{eqn:aspinNewtonStep}.
The vector $g_i^{\text{ASPIN}}$ is a substitution for the right-hand side~$g_i$ of \eref{eqn:newtonStep}.
If we chose $g_i^{\text{ASPIN}} = C^{-1}_i g_i$, we would obtain a left preconditioned system. 
However, in the ASPIN method, $g_i^{\text{ASPIN}}$ is defined differently.
It is the sum of local contributions $\locals \in \Sk$ prolongated to $\mathbb R^n$ via the prolongation operator $\Ik$:
\begin{equation}
\label{eqn:preconditionedGradientASPIN}
	g_i^{\text{ASPIN}} = -\sum_k \Ik \locals.
\end{equation}
Here, each local contribution $\locals \in \Sk$ is the solution of a non-linear system of equations
\begin{equation}
\label{eqn:aspinStep}
	\nabla \Hk(\Pk \globiti+\locals) = 0 \qquad \forall k\in\{1,\ldots,N\},
\end{equation}
where $\Pk$ is the projection operator defined in \eref{eqn:defPk} and $\Hk$ is a local objective function.
Here, instead of restricting the  original objective function to the subdomain $\Sk$ in order to obtain $\Hk$ as done in the original work~\cite{CaiKeyes00},
we advocate to employ the following, general objective function, cf.~\cite{GrattonSartenaerToint06,Nash00}, 
\begin{equation}
	\label{eqn:exactLowerLevelModel2}
	\Hk(\Pk\globiti+\locals) = 
		\Jk(\Pk\globiti+\locals) + 
		\langle \delta g^k, \locals \rangle,	
\end{equation}	
where $\langle \cdot,\cdot \rangle$ is the Euclidean inner product, $\Jk$ is some arbitrary, sufficiently smooth local objective function and
$\delta g^k = \Rk g_i - \nabla \Jk(\Pk\globiti)$ denotes the difference between the restricted global
gradient~$\Rk g_i$ and the initial subset gradient~$\nabla \Jk(\Pk\globiti)$.
Here, as opposed for example to \cite{GrattonSartenaerToint06}, we use different restriction operators for $u_i$ and $g_i$.
Once a choice for the local objective functions~$\Jk$ is made, 
we thus have defined the right-hand side $g_i^{\text{ASPIN}}$ of \eref{eqn:aspinNewtonStep}.

An important feature of the general objective function~$\Hk$ is that it is first-order consistent with 
the global objective function~$J$ regardless of the definition of the global objective functions~$\Jk$. 
This means that the following asymptotic result holds
\begin{equation*}
	\nabla \Hk(\localu+\locals) \to \Rk g_i \qquad \text{ if } \qquad \Pk\globiti+\locals \to \Pk \globiti.
\end{equation*}

It remains to define the global preconditioning matrix~$C_i^{-1}$, which appears on the left-hand side of \eref{eqn:aspinNewtonStep}.
Following \cite{CaiKeyes00},
in the ASPIN framework the  matrix~$C_i^{-1}$ is defined as (a general variant) of the
additive Schwarz preconditioner
\begin{equation}
\label{eqn:matrixC}
	C_i^{-1} = \sum_k (\Ck)^{-1} = \sum_k \left[\Ik \left(\Bk(\Pk \globiti) \right)^{-1} \Rk\right],
\end{equation}
where $\Bk(\Pk \globiti)$ are invertible approximations to $\nabla^2 \Hk(\Pk \globiti)$.

\section{A Short Survey of Trust-Region Methods}
\label{sec:trustRegion}
While the Newton and ASPIN methods presented in Section~\ref{sec:aspin} can provide fast convergence to a first-order critical point in the vicinity of a solution under certain assumptions,
for an arbitrary initial iterate~$\globitstart$ the convergence of the method cannot be guaranteed, cf.~\cite[Section~3.3]{NocedalWright2006}.
In this section, we therefore briefly present the trust-region framework, cf.~\cite{ConnGouldToint00},
which will be coupled with the ASPIN approach in oder to obtain a solution approach that converges for arbitrary starting iterates
and reduces to the fast ASPIN method in the vicinity of a solution.

Trust-region methods, also known as \emph{restricted step} methods, are used to solve optimization problems of the kind~\eqref{eqn:nlMinProb} employing the following strategy.
Within a subset of the region, the \emph{trust region}, the objective function~$J$ is approximated using a model function~$\psi$, often a quadratic one.
A new iterate is computed by solving an optimization problem using this model,
which is supposed to be easier to solve than the original optimization problem.
If the model describes the reduction in the objective function~$J$ within the trust region sufficiently accurately, 
then the new iterate is accepted and the region is expanded, 
otherwise it is rejected and the trust region is contracted.

More precisely, for our setting, for a given iterate~$\globiti \in \mathbb{R}^n$ the \emph{trust-region model}~$m_i$
is a quadratic approximation to the function~$J$ at the point~$\globiti$. 
It is given as 
$m_{i}: \mathbb R^n \to \mathbb R$ with
\begin{equation*}
	m_i(\globiti + s) = m_i(\globiti) + \psi_{i}(s) \quad \text{for } s \in \mathbb{R}^n,
\end{equation*}
where $m_i(\globiti) = J(\globiti)$ and
$\psi_{i}: \mathbb R^n \to \mathbb R$	
is defined as
\begin{equation}
	\label{eqn:localQuadraticModel}
	\psi_{i}(s) = \langle g_{i}, s\rangle+\frac{1}{2}\langle s, B_{i} s \rangle \quad \text{for } s\in \mathbb{R}^n.
\end{equation}
Here, $g_i = \nabla J(\globiti)$ is the gradient and $B_i=B(\globiti)\in \mathbb{R}^{n\times n}$ is a symmetric approximation to the Hessian matrix~$ \nabla^2 J(\globiti)$.
We then have that for given $s \in \mathbb{R}^n$, the reduction in the model, the \emph{predicted reduction}, is
\begin{equation*}		
	\pred_{i}(s) = m_i(\globiti) - m_i(\globiti + s) =  -\psi_{i}(s).
\end{equation*}	
The function $-\psi_i(s)$ is a first-order -- if we choose $B_i = \nabla^2 J(\globiti)$, a second-order -- Taylor approximation to 
the \emph{actual reduction}
\begin{equation*}		
	\ared_{i}(s) = J(u_{i})-J(u_{i}+s).
\end{equation*}

Then, a trust-region correction $s_{i} \in \mathbb R^n$ minimizes the model in the trust region,
and is defined as the solution 
of the following constrained, quadratic programming problem, called \emph{trust-region subproblem}
\begin{equation}
	\label{eqn:localMinProb}
	\globstepi := \argmin_{s \in \mathbb{R}^n,\, \norm{s} \le \Delta_{i}^G} m_i(\globiti + s) = \argmin_{s \in \mathbb{R}^n,\,\norm{s} \le \Delta_{i}^G} \psi_{i}(s),
\end{equation}
where $\Delta_{i}^G \in \mathbb R^+$ is called the \emph{trust-region radius}. 
Problem~\eqref{eqn:localMinProb} is not necessarily solved accurately, that is,
$\globstepi$ is not necessarily an exact solution of \eref{eqn:localMinProb}. 
It can rather be an approximation to the minimizer, 
as long as $\globstepi$ satisfies the following \emph{sufficient decrease condition}
\begin{equation}
\label{eqn:sufficientDecrease}
-\psi_i(\globstepi)\ge \beta_1 \min\left\{\|g_i\|^2, \|g_i\| \Delta_i^G \right\},
\end{equation}
where $\beta_1>0$ is an appropriately chosen constant. 
The two obvious choices for the step $\globstepi$ are the Cauchy point, that is 
the point minimizing the quadratic model in the steepest descent direction subject to the step being within the trust region,
and the point that makes the model $\psi_i$ as small as possible within the trust region. 
In practice, any point lying between these two extremes should be acceptable~\cite{ConnGouldToint00}.
Throughout this paper we assume that $\beta_1$ is chosen in such a way that also the worst of these possibilities, the Cauchy point, satisfies Condition~\eqref{eqn:sufficientDecrease},
and that the correction~$\globstepi$ is at least as good as the Cauchy point.

\begin{remark}
\label{rem:relationMinProbSystem}
Note that the solution of the trust-region subproblem is closely related to the Newton step~\eqref{eqn:newtonStep}.
The vector $\globstepi$ is a global solution of the trust-region subproblem~\eqref{eqn:localMinProb},
if and only if $\globstepi$ is feasible and there is a scalar $\lambda \geq 0$ such that the following conditions are satisfied~\cite[Theorem~4.1]{NocedalWright2006}:
\begin{align*}
	\left(B_i + \lambda I\right) s_i &= - g_i \\
	\lambda \left( \Delta_{i}^G - \norm{\globstepi}\right) &=0\\
	\left(B_i + \lambda I\right)& \text{ is positive semidefinite.}
\end{align*}
\end{remark} 

In a second step, the \emph{decrease ratio}~$\rho_{i}$ is computed
comparing the predicted with the actual reduction.
It is defined as
\begin{equation}	
\label{eqn:decreaseRatio}
	\rho_{i} = 
		\frac{\ared_{i}(\globstepi)}{\pred_{i}(\globstepi)},
\end{equation}
and used to update the trust-region radius~$\Delta_i^G$ 
and to accept or reject the correction. 
That is, one defines
\begin{equation}
\label{eqn:acceptanceTrustRegion}
	u_{i+1} = \begin{cases} 
			\globiti+\globstepi & \text{if } \rho_{i} \ge \eta,\\
			\globiti	& \text{otherwise,}
	          \end{cases}
\end{equation} 
and
\begin{equation}
	\label{eqn:updateTrustRegionRadius}
	\Delta_{i+1}^G = 
		\begin{cases}
			\gamma_2\Delta_{i}^G & 
				\text{if\ }  \rho_{i} \ge \eta,\\
			\gamma_1 \Delta_{i}^G & 
				 \text{otherwise,}
		\end{cases}
\end{equation}
where $1 > \eta > 0$ and $\gamma_2 > 1 > \gamma_1 > 0$ are constants
assumed to be given a priori and fixed for the whole computation. 
These four steps are summed up in Algorithm~\ref{alg:TrustRegionSolver}.

The convergence analysis of trust-region methods can be carried out based
on the following moderate assumptions, for details see~\cite{ColemanLi94,ConnGouldToint00}.
We state these assumptions here, as we will also base our following convergence analysis on them.
\begin{itemize}
	\item[(\ATR1)]{
		For a given initial iterate $\globitstart \in \mathbb R^n$, %in Algorithm~\ref{alg:gaspin},
		we assume that the level set
		\begin{equation*}
		\mathcal L = \{u\in \mathbb R^n \mid J(u) \le J(\globitstart)\}
		\end{equation*}
		is compact.
	}
	\item[(\ATR2)]{
		We assume that $J$ is continuously differentiable on $\mathcal L$, 
		and that the norms of the gradients are bounded 
		by a constant $C_g>0$, i.e., $\|\nabla J(u)\| \le C_g$ for all $u\in \mathcal L$.
	}

	\item[(\ATR3)]{
		There exists a constant
		$C_B > 0$ such that for all iterates $u \in \mathcal L$ and for each symmetric matrix
		$B(u)$ employed in \eref{eqn:localQuadraticModel} the inequality
		$\|B(u)\| \le C_B$ is satisfied.
	}
\end{itemize}

\begin{algorithm}
	\caption{Trust-Region Algorithm}
		\label{alg:TrustRegionSolver}
	\begin{algorithmic}[1]
		\Statex
		\Input{
			$J:\mathbb R^n \to \mathbb R$, $\globitstart \in \mathbb R^n$, $\Delta_0^G \in \mathbb R^+$,
			$\gamma_1$, $\eta\in (0,1)$, $\gamma_2>1$, $m\in \mathbb N\cup\{\infty\}$
		}
		\Output{
			Iterate~$u_m$
		}
		\Statex
			\For {$i \gets 0,1,\ldots, m-1$}
				\State Generate $\psi_i$ using \eref{eqn:localQuadraticModel}.
				\State Solve Problem~\eqref{eqn:localMinProb} approximately such that  \eqref{eqn:sufficientDecrease} holds,
				and obtain $\globstepi \in \mathbb R^n$.
				\State Compute $\rho_i$ according to \eref{eqn:decreaseRatio}.
				\If {$\rho_i \ge \eta$}
       				\State $u_{i+1} \gets \globiti+\globstepi$
       			\Else
       				\State $u_{i+1} \gets \globiti$
    			\EndIf
    			\State Compute $\Delta_{i+1}^G$ according to \eref{eqn:updateTrustRegionRadius}.
			\EndFor
			\State \Return $u_{m}$
	\end{algorithmic}
	
\end{algorithm}

\input{gaspin_definition.tex}

\section{Convergence of the Nonlinearly Left Preconditioned Trust-Region Strategy}
\label{sec:firstOrderConvergence}

In the present section, we analyze the convergence properties of Algorithm~\ref{alg:gaspin}. 
In particular, 
we will show that this algorithm generates a sequence of iterates converging to first-order critical points 
under the same assumptions as used in Section~\ref{sec:trustRegion}. 

Note that in contrast to \cite{Toint1988}, we do not assume that $\|\tilde g_i-g_i\|\to 0$ but employ the 
modified sufficient decrease condition \eref{eqn:modifiedSuffDecrease}.
Furthermore, in contrast to \cite{CGrossRKrause2009a},
we do not state further assumptions on the local objective functions $\Hk$ and the local solution 
process. As we will see, this is not necessary since Assumption~\eqref{eqn:localConstraint} is sufficiently
strong to ensure convergence to critical points. 
On the other hand, let us remark that the proof of 
Theorem~\ref{thm:limInfTr} will be carried out by contradiction, i.e., by assuming that
$\|g_i\| \ge \varepsilon >0$. Therefore, in the following lemma we show that in this case, for 
sufficiently small $\Delta_i^L$ or $\Delta_i^G$, Assumption~\eqref{eqn:modifiedSuffDecrease2} holds.

\begin{lem}
\label{lem:sufficientDecrease}
	Assume that assumptions (\ATR1), (\ATR2), (\ATR3) hold. 
	Furthermore, assume that there exists an 
	$\varepsilon>0$ such that $\|g_i\|\ge\varepsilon$ and that either $\Delta_i^L$ or $\Delta_i^G$ is sufficiently
	small. Then every correction $\globstepi$ computed in Algorithm~\ref{alg:gaspin} satisfies 
	\eref{eqn:modifiedSuffDecrease2}.
\end{lem}
\begin{proof}
	Throughout the proof of this lemma, we assume that $\Delta_i^L$ is sufficiently small.
	If instead $\Delta_i^G$ is sufficiently small, $\Delta_{i}^L$ will also be sufficiently small by \eqref{eqn:boundDeltaL} and \eqref{eqn:deltaGUpdate}. 

	Due to Assumption~\eqref{eqn:assumptionDeltaS}, we have that 	
	\begin{equation*}
		\Delta_i^L \ge \|g_i-\tilde g_i\| \ge \|g_i\|-\|\tilde g_i\|.
	\end{equation*}	
	Thus, we obtain $\|\tilde g_i\| \ge \|g_i\|-\Delta_i^L$ and if we assume that 
	$\Delta_i^L$ is sufficiently small, we obtain
	\begin{equation}
		\label{eqn:boundBelowgiTilde}
		\|\tilde g_i\| \ge \|g_i\|-\Delta_i^L \ge 0.
	\end{equation}			
	Now we investigate \eref{eqn:modifiedSuffDecrease2}, i.e., 
	\begin{equation*}
		\beta_1 \min\{\|\tilde g_i\|^2, \|\tilde g_i\| \Delta_i^G \} 
		\ge \beta_2 \min\{\|g_i\|^2, \|g_i\| \Delta_i^G \},		
	\end{equation*}
	and determine $\Delta_i^L$ controlling the difference between $g_i$ and $\tilde g_i$ such that this inequality holds.
	To this end, we make the following case differentiation:
	\begin{description}
	\item[1. Assume that $\norm{\tilde g_i} \leq \Delta_i^G$ and $\norm{g_i}\leq \Delta_i^G$.]
		Then \eref{eqn:modifiedSuffDecrease2} is equivalent to
		\begin{equation*}
			\beta_1 \norm{\tilde g_i}^2
			\ge \beta_2 \norm{g_i}^2.
		\end{equation*}
		Using \eref{eqn:boundBelowgiTilde}, in order for \eref{eqn:modifiedSuffDecrease2} to hold, it is thus sufficient to find a parameter $\Delta_i^L$, such that
		\begin{equation*}
			\beta_1 \left(\norm{ g_i} -\Delta_i^L \right)^2
			\ge \beta_2 \norm{g_i}^2
		\end{equation*}
		holds.
		In order to determine such $\Delta_i^L$, we solve 
		$$
		\Delta \in \mathbb R: \beta_2\|g_i\|^2 = \beta_1(\|g_i\|-\Delta)^2
		$$ 
		and obtain
		\begin{equation*}
			\Delta = \left(1 \pm \sqrt{\frac{\beta_2}{\beta_1}}\right)\|g_i\|.
		\end{equation*}
		Thus \eref{eqn:modifiedSuffDecrease2} holds  in this case if 
		$0< \Delta_i^L \le \left(1-\sqrt{\frac{\beta_2}{\beta_1}}\right)\norm{g_i}$.
	\item[2. Assume that $\norm{\tilde g_i} \geq \Delta_i^G$ and $\norm{g_i}\geq \Delta_i^G$.]
		Then \eref{eqn:modifiedSuffDecrease2} is equivalent to
		\begin{equation*}
			\beta_1 \norm{\tilde g_i} \Delta_i^G
			\ge \beta_2 \norm{g_i} \Delta_i^G.
		\end{equation*}
		Using \eref{eqn:boundBelowgiTilde}, in order for \eref{eqn:modifiedSuffDecrease2} to hold, it is thus sufficient to find a parameter $\Delta_i^L$, such that
		\begin{equation*}
			\beta_1 \left(\norm{ g_i} -\Delta_i^L \right)
			\ge \beta_2 \norm{g_i} 
		\end{equation*}
		holds, where we additionally divided both sides of the inequality by $\Delta_i^G>0$.
		Consequently, \eref{eqn:modifiedSuffDecrease2} holds if 
		\begin{equation*}
			0<\Delta_i^L \le \left(1-\frac{\beta_2}{\beta_1}\right)\|g_i\|,
		\end{equation*}
		where $\left(1-\frac{\beta_2}{\beta_1}\right)\in (0,1)$.
	\item[3a. Assume that $\norm{\tilde g_i} < \Delta_i^G$ and $\norm{g_i} > \Delta_i^G$.]
		We remark that this state can be considered as an intermediate state, since for $\Delta_i^G\to 0$ or 
		just for $\Delta_i^L\to 0$ we obtain $\tilde g_i\to g_i$ where $\|g_i\|\ge\varepsilon$ and
		eventually 
		\begin{equation*}
			\Delta_i^G \leq \norm{\tilde g_i}
		\end{equation*}
		and we are in Case~2.
		Note that this reasoning is feasible, since $u_i$ is only updated, and thus $g_i$ is only changed, if \eref{eqn:modifiedSuffDecrease2} holds, 
		cf.~\eqref{eqn:updateU}.
		
	\item[3b. Assume that $\norm{\tilde g_i} > \Delta_i^G$ and $\norm{g_i} < \Delta_i^G$.]
	 	This is the second intermediate state.
		Here, as well, we have for $\Delta_i^L\to 0$ that $\tilde g_i\to g_i$ and thus that Case~2 eventually holds,
		if $\Delta_i^G$ is fixed.
		If $\Delta_i^G \to 0$, we will eventually reach either Case~3a or Case~1. 
		
	\end{description}	 
	As Cases~3a and 3b show, either Case~1 or Case~2 holds if $\Delta_i^L$ or $\Delta_i^G\ge\Delta_i^L$ become sufficiently small.
	Then, we have that there exists some $\Delta_i^L>0$ (independent from $\Delta_i^G$)
	which satisfies \eref{eqn:modifiedSuffDecrease2}. 
	This proves the proposition.

\end{proof}

\begin{remark}
	We remark that in the intermediate cases 3a and b of the previous proof, 
	we are not able to compute $\Delta_i^L$ employing the same trick as for the first two cases.
	In Case 3a, by solving 
	$$
		\Delta \in \mathbb R: \beta_2\|g_i\|\Delta_i^G  = \beta_1(\|g_i\|-\Delta)^2,
	$$
	we obtain that 
	\begin{equation*}
		\|g_i\|-\left(\|g_i\|^2+\frac{\beta_2}{\beta_1}\|g_i\|\Delta_i^G\right)^{\frac 1 2} \le \Delta_i^L \le
		\|g_i\|+\left(\|g_i\|^2+\frac{\beta_2}{\beta_1}\|g_i\|\Delta_i^G\right)^{\frac 1 2}.
	\end{equation*}
	This contradicts $\Delta_i^L\le \Delta_i^G$ from \eref{eqn:boundDeltaL} and $\|g_i\|-\Delta_i^G \ge 0$.
	
	In Case 3b, by solving 
	$$
		\Delta \in \mathbb R: \beta_1(\|g_i\|-\Delta)\Delta_i^G  = \beta_2\|g_i\|^2,
	$$ 
	we obtain that 
	\begin{equation*} 
		 -\frac{\beta_2}{\beta_1 \Delta_i^G}\|g_i\|^2+\|g_i\|\ge \Delta_i^L,
	\end{equation*}
	which in general contradicts $\Delta_i^L \in \mathbb R^+$.
\end{remark}

The previous lemma shows that a sufficient decrease in the objective function is possible, if the perturbation
in the quadratic model becomes small enough. Now, we exploit a quite similar argumentation and 
the mean value theorem to show that the quadratic approximation~$\tilde \psi_i$ to the actual decrease
becomes asymptotically exact. 
% Therefore, in the following lemma we prove that the predicted reduction becomes sufficiently accurate 
% if $\Delta_i^G$ becomes sufficiently small.

\begin{lem}
\label{lem:rhoToOne}
	Let assumptions (\ATR1), (\ATR2) and (\ATR3) hold. Suppose, moreover, that there exists an $\varepsilon > 0$ such that 
	$\| g_{i}\| \ge \varepsilon > 0$ and that $\Delta^G_{i}$ is sufficiently small.
	Then the decrease ratio defined by \eref{eqn:decreaseratio} satisfies
	\begin{equation*}
			\tilde \rho_{i} \ge \eta.
	\end{equation*}		
\end{lem}

\begin{proof}	
	Exploiting (\ATR1), (\ATR2) and the mean value theorem yields
	\begin{equation*}
		J(u_{i}+s_{i})-J(u_{i}) = \langle \overline g_{i}, s_{i}\rangle
	\end{equation*}
	with $\overline g_{i} = \nabla J(u_{i} + \tau s_{i})$ with $\tau\in(0,1)$. 		
	Using the definitions of the decrease ratio and $\widetilde \psi_{i}$, as well as 
	(\ATR2) and (\ATR3) yields
	\begin{equation}
	\label{eqn:estimationRho}
	\begin{aligned}
	|\widetilde{\psi}_i(s_{i})||\tilde \rho_i-1|
			=&\abs{J(u_{i}+s_{i}) - J(u_{i})- \langle \tilde g_{i}, s_{i}\rangle - \tfrac{1}{2}\langle s_{i}, B_{i} s_{i}\rangle}\\
			\le& \tfrac{1}{2} \abs{\langle s_{i}, B_{i} s_{i} \rangle } + 
					\abs{\langle \overline g_{i}-\tilde g_{i}, s_{i}\rangle}\\
% 			\le \frac{1}{2} C_B \|s_{i}\|^2+\|\overline g_{i}-\tilde g_{i}\|\|s_{i}\| \\[0.2cm]
			\le& \tfrac{1}{2} C_B \|s_{i}\|^2+\|\overline g_{i}-\tilde g_{i}\| \|s_{i}\|\\
			\le&  \tfrac{1}{2} C_B (\Delta^G_i)^2+\|\overline g_{i}-\tilde g_{i}\| \Delta^G_i.
	\end{aligned}
	\end{equation}
	Due to \eref{eqn:assumptionDeltaS} and \eref{eqn:boundDeltaL}, we have in particular 
	\begin{equation}
	\label{eqn:estimationMean}
		\norm{\overline g_{i}-\tilde g_{i}} 
		= \norm{\overline g_i- g_i+ g_i - \tilde g_i } 
		\le \|\overline g_{i}-g_{i}\|+\Delta_i^L
		\le \|\overline g_{i}-g_{i}\|+\Delta_i^G.
	\end{equation}	
	Following Lemma~\ref{lem:sufficientDecrease}, if $\Delta^G_i$ and thus $\Delta_i^L$ are sufficiently small, \eref{eqn:modifiedSuffDecrease} yields
	\begin{equation}	
	\label{eqn:suffDecInRhoEst}
		-\widetilde{\psi}_i(\globstepi) \ge \beta_1 \min\left\{\|\tilde g_i\|^2, \|\tilde g_i\| \Delta_i^G \right\} \ge \beta_2 \min\left\{\|g_i\|^2, \|g_i\| \Delta_i^G \right\} > 0.
	\end{equation}	
	
	Combining the previous results, we now show that $\abs{\tilde \rho_i -1} \to 0$ for $\Delta_i^G \to 0$.
	With $\varepsilon >0$ from the assumptions of the lemma, we obtain that
	\begin{align*}
		\beta_2\varepsilon|\tilde\rho_i-1|
		=& (\Delta_{i}^G )^{-1} \beta_2 \varepsilon \min \left\{\varepsilon, \Delta^G_i\right\} \abs{\tilde \rho_{i} -1} &\text{if $\Delta_i^G$ is sufficiently small}\\
		\le& (\Delta_{i}^G )^{-1} \beta_2 \min \left\{\norm{g_i}^2, \norm{g_i}\Delta^G_i \right\} \abs{\tilde \rho_{i} -1}&\text{as $\norm{g_i}\geq\varepsilon$ by assumption}\\
		\le& (\Delta^G_i )^{-1} \abs{\tilde \psi_{i}(\globstepi)} \abs{\tilde \rho_{i} -1} &\text{by \eref{eqn:suffDecInRhoEst}}\\
		\le &	\frac{ C_B}{2}\Delta^G_i  + \|\overline g_{i} - \tilde g_{i} \|  &\text{by \eref{eqn:estimationRho}}\\
		\le &	\frac{ C_B}{2}\Delta^G_i  + \|\overline g_{i} - g_{i} \| + \Delta^G_i &\text{by \eref{eqn:estimationMean}.}
	\end{align*}
	We conclude that if $\Delta^G_i$ converges to zero,
	the right-hand side of this inequality goes to zero as well. This is due to the fact that for 
	$\Delta^G_i\to0$ also 
	$\|s_{i}\|\to 0$ and $(u_{i})_i$ converges in $\mathcal L$. 
	Thus, $\abs{\tilde \rho_i -1}$ is bounded from above by a term that converges to zero for $\Delta_i^G \to 0$.
	Therefore, if $\Delta^G_i$ is sufficiently small, we obtain that
	\begin{equation*}
		\tilde \rho_i \ge \eta
	\end{equation*}
	holds.
\end{proof}

One observation in the proof of the previous lemma is that the perturbation of the gradient can be estimated by $\mathcal O(\Delta_i^G)$, 
cf.\ for instance~\cite{ConnEtAl1993}. 
We conclude from Lemma~\ref{lem:rhoToOne}, that the smaller the trust-region radius~$\Delta_i^G$ gets, the more accurate the preconditioned model becomes. 
This is the final step for proving convergence of our globalized ASPIN method, Algorithm~\ref{alg:gaspin} in the following theorem. 

\begin{thm}
\label{thm:limInfTr}
	Let assumptions (\ATR1), (\ATR2) and (\ATR3) hold. 
	In this case, we obtain that the sequence of iterates generated by Algorithm \ref{alg:gaspin} has the property
	\begin{equation*}
		\liminf_{i\to\infty} \| g_i\| = 0.
	\end{equation*}
\end{thm}
\begin{proof}
	Assume that the proposition does not hold, i.e., there exists an $\varepsilon>0$ and an index 
	$\nu_0$ such that $\| g_i\| > \varepsilon$ for all $i\ge\nu_0$. 
	We will show that if this is the case, the sequence of trust-region radii converges to zero.

	If there are only finitely many successful corrections, 
	the update criteria \eqref{eqn:deltaGUpdate} and \eqref{eqn:boundDeltaL} directly imply that
	$\Delta^L_i\to 0$. 
	Then also $\Delta_i^G \to 0$, since the case that only $\Delta_i^L$ is reduced in each iteration and $\Delta_i^G$ stays constant does not take place.
	In fact, if Assumption~\eqref{eqn:modifiedSuffDecrease2} subsequently does not
	hold we obtain $\Delta^L_i \to 0$. Then Lemma~\ref{lem:sufficientDecrease} 
	gives that \eref{eqn:modifiedSuffDecrease2} will (after finitely many steps) hold, too.
	Since we have only finitely many successful iterations, we obtain $\Delta^G_i \to 0$.

	If there are infinitely many successful corrections, we have for every successful correction due to \eref{eqn:modifiedSuffDecrease} and $\tilde \rho_i > \eta$ that
	\begin{equation*}
		J(\globiti)-J(u_{i+1}) \ge \eta \beta_2 \varepsilon \min \{\varepsilon, \Delta_i^G\}.
	\end{equation*}
	Since by (\ATR1) the levelset $\mathcal L$ is compact, the sequence 
	$\left(J(u_k)\right)_k$ is non-increasing and bounded from below, and thus a Cauchy sequence.
	We obtain consequently
	\begin{equation*}
		J(\globiti)-J(u_{i+1}) \to 0,
	\end{equation*}
	which implies $\Delta^G_i \to 0$.

	Then we can use Lemmas~\ref{lem:sufficientDecrease} and~\ref{lem:rhoToOne} and obtain that eventually, for sufficiently small $\Delta^L_i$
	and $\Delta^G_i$, every correction is successful.
	This contradicts $\Delta^L_i, \Delta^G_i \to 0$ and proves the proposition.
\end{proof}

\begin{thm}
\label{thm:firstOrderCriticalTr}
	Let assumptions (\ATR1), (\ATR2) and (\ATR3) hold. 
	Then the sequence of iterates generated by Algorithm \ref{alg:gaspin} 
	converges to a first-order critical point, i.e.,
	\begin{equation*}
		\lim_{i\to\infty} \| g_i\| = 0.
	\end{equation*}
\end{thm}
\begin{proof}
	This proof follows exactly the same reasoning as the proof of \cite[Theorem~6.6]{UlbrichUlbrichHeinkenschloss99}.
	
\end{proof}

\section{Numerical Examples}
\label{sec:numerics}

\newcommand{\x}{\boldsymbol{x}}
\renewcommand{\u}{\boldsymbol{u}}
\newcommand{\p}{\boldsymbol{p}}
\newcommand{\f}{\boldsymbol{f}}
\newcommand{\F}{\boldsymbol{F}}
\newcommand{\G}{\boldsymbol{G}}
\newcommand{\C}{\boldsymbol{C}}
\newcommand{\n}{\boldsymbol{n}}
\renewcommand{\div}{\text{div}}
\newcommand{\T}{\hat{\boldsymbol{T}}}
\newcommand{\vp}{\boldsymbol{\varphi}}
\newcommand{\W}{\tilde{\boldsymbol{W}}}
\newcommand{\nvp}{\nabla \boldsymbol{\varphi}}
\newcommand{\cof}{ \text{\emph{Cof}}\/}
\newcommand{\idet}{ \text{\emph{det}}}
\newcommand{\E}{\boldsymbol{E}}

We employ the globalized ASPIN strategy for the solution of 
optimization problems arising from the field of non-linear elasticity. In our applications,
we are interested in the computation of energy optimal displacements, where we 
follow~\cite{Ciarlet88} and employ the following polyconvex stored energy function
\begin{equation}
\label{eqn:ogdenFuncApts}
W(\x,\C)= 3(a+b)+(2a+4b) \cdot \text{tr} \E+2b\cdot(\text{tr} \E)^2-2b\cdot\text{tr}(\E^2) + \Gamma(\text{det} (\nvp))	
\end{equation}
for $\x\in\Omega\subset \mathbb R^d$, $d=2,3$. 
Here, $\C=\C(\u) = (I+\nabla \u)^T(I+\nabla \u)$ is the right Cauchy-Green strain tensor,
$\E=\E(\u)=\frac{1}{2}(\C(\u)-I)$ is the Green-St.\ Venant strain tensor, $\nabla\vp=\text{Id}+\nabla\u$ is the deformation tensor and
$\Gamma(\delta)=c\delta^2-d\log \delta$ is a logarithmic barrier function.
For our examples, the constants are chosen as follows:
\begin{equation}
\label{eqn:constants}
a = \mu+\frac{1}{2}\Gamma'(1) \text{, } b = -\frac{\mu}{2}-\frac{1}{2}\Gamma'(1) \text{, } c = -\frac{\lambda}{4}-\mu \text{ and } d=\frac{3\lambda}{4}+\mu,
\end{equation}	
where $\lambda$ and $\mu$ are the Lam\'e constants.

A particular and important property of this class of stored energy functions is that (depending on the choice of $B_i$) assumptions
(\ATR1)-(\ATR3) hold, cf.~\cite{CGross_2009}. 
Therefore, by Theorem~\ref{thm:firstOrderCriticalTr},
the globalized ASPIN strategy, Algorithm~\ref{alg:gaspin}, 
 provably computes a first-order critical point of 
\begin{equation*}
	\u\in S_h: J(\u) = \int_\Omega W(\x,\C) + f\cdot \u d\x + \int_{\Gamma_N} \p \cdot \u da = \min!
\end{equation*}
for given Dirichlet values\footnote{Prescribed displacements at some boundaries} at $\Gamma_D \subseteq \partial \Omega$, Neumann boundary conditions 
on $\p\in [C(\Gamma_N)]^3$, where $\Gamma_N = (\partial\Omega\setminus \Gamma_D)$, and volume forces $\f \in [C(\Omega)]^3$. Here, the space of linear
finite elements is denoted by $S_h$.

\subsection{Implementational Aspects and Runtime Comparisons}

The algorithmic framework presented in this article was implemented in {\sc ObsLib++}, a framework for the solution of 
constrained optimization problems arising from the finite element discretization of elastic PDEs~\cite{Obslib_2006,Gross08}. 
{\sc ObsLib++} employs, as a grid manager, the parallelized unstructured grid manager UG \cite{PBastian_KBirken_KJohannsen_1997}, 
which was extended in order to allow for asynchronously applied  trust-region and linesearch methods \cite{CGross_2009}.

\subsubsection{Evaluating the Extended Sufficient Decrease Condition}

In our implementation, the evaluation of the modified sufficient decrease condition~\eqref{eqn:modifiedSuffDecrease}  is a two-step process. 
First, we compute the Cauchy point~$\tilde s_i^C$ 
for $\widetilde \psi_i$. Following the argumentation in \cite{ConnGouldToint00}, if (\ATR1), (\ATR3) hold for a quadratic model and if $\tilde g_i\not=0$, 
a Cauchy point induces a sufficient decrease of a quadratic model
with a constant $\tilde \beta_1 = \tilde \beta_1(C_B)>0$, that is
\begin{equation*}
	-\widetilde \psi_i\left(\tilde s_i^C\right) \geq \tilde \beta_1 \min \left\{ \norm{\tilde g_i}^2, \norm{\tilde{g_i}}\Delta_i^G \right\}.
\end{equation*}
Therefore, if the correction vector $\globstepi$ satisfies
\begin{equation*}
	-\widetilde \psi_i(\globstepi) \ge -c_1 \widetilde \psi_i(\tilde s_i^C) 
\end{equation*}
for $c_1>0$, 
% e.g., $c_i=10^{-4}$, 
we know that $\globstepi$ satisfies  \eref{eqn:modifiedSuffDecrease1} with $\beta_1=c_1 \tilde \beta_1$.

In the second step, we compute the Cauchy point~$\globstepi^C$ for $\psi_i$. Here, we have that $\globstepi^C$ satisfies the the second equation in the sufficient decrease condition~\eqref{eqn:sufficientDecrease}
with a constant $\tilde \beta_2 = \tilde \beta_2(C_B)>0$, that is
\begin{equation*}
	-\psi_i(\globstepi^C) \geq \tilde \beta_2 \min \left\{ \norm{ g_i}^2, \norm{{g_i}}\Delta_i^G \right\}.
\end{equation*}
Therefore, if $\tilde g_i \not=0$ and if 
\begin{equation}
\label{eqn:cauchyPoints}
	-c_1 \widetilde \psi_i(\tilde \globstepi^C) \ge -c_2 \psi_i(\globstepi^C)
\end{equation}
for $c_2>0$,
$\tilde \globstepi^C$ satisfies \eref{eqn:modifiedSuffDecrease} 
with the constant $\beta_2 = c_2 \tilde \beta_2$.

\begin{figure}
\centering
\includegraphics[width=0.42\textwidth]{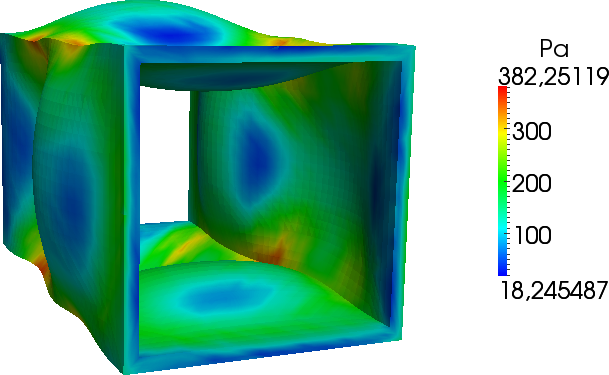}\hspace{0.05\textwidth}
\includegraphics[width=0.42\textwidth]{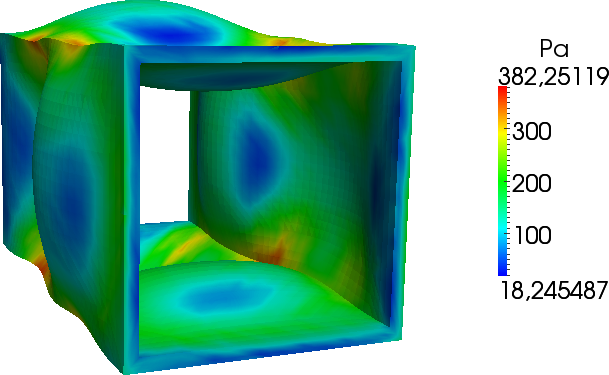}\\[5mm]
\includegraphics[width=0.42\textwidth]{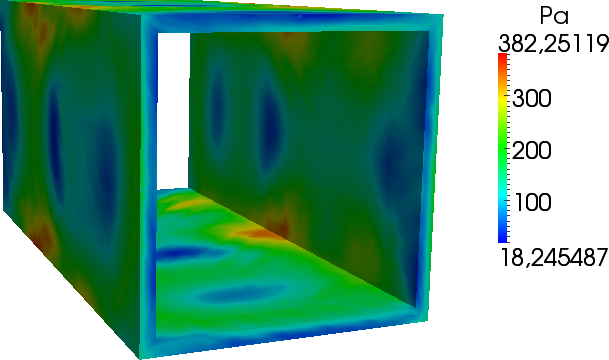}

\caption{\label{fig:brick}Numerical results of the example from Section \ref{sec:brick} with 100,800 unknowns, where
the initial geometry is deformed according to the solution of \eref{eqn:nlMinProb} (\emph{upper figures}). \emph{Lower
figure}: the initial geometry. Colors denote the von-Mises stresses.
}
\end{figure}

\begin{figure}
	\centering
	{\bf Computation on 240 Processors: Trust-Region vs. Globalized ASPIN}\\
	\includegraphics[width=0.45\textwidth]{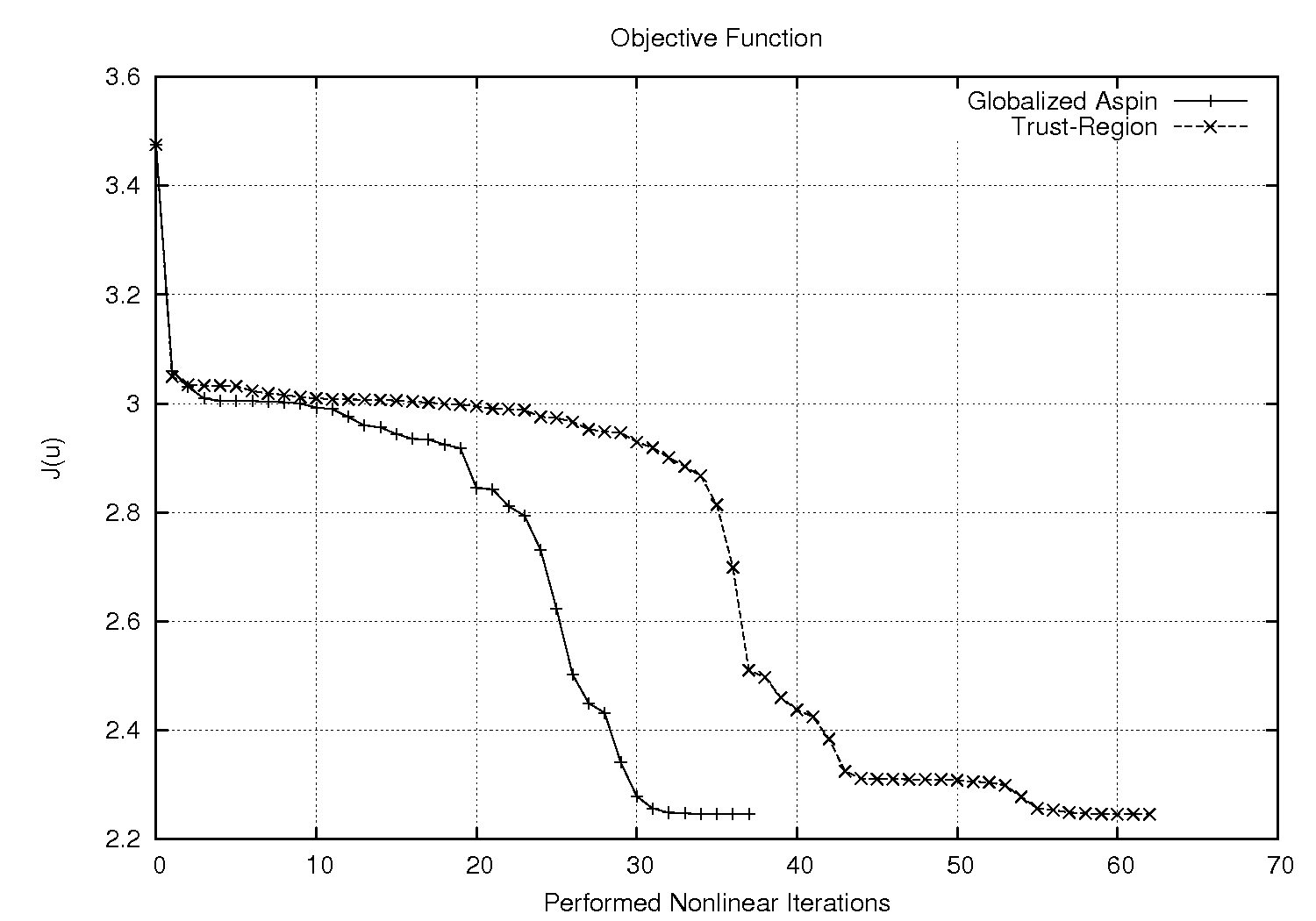}
	\includegraphics[width=0.45\textwidth]{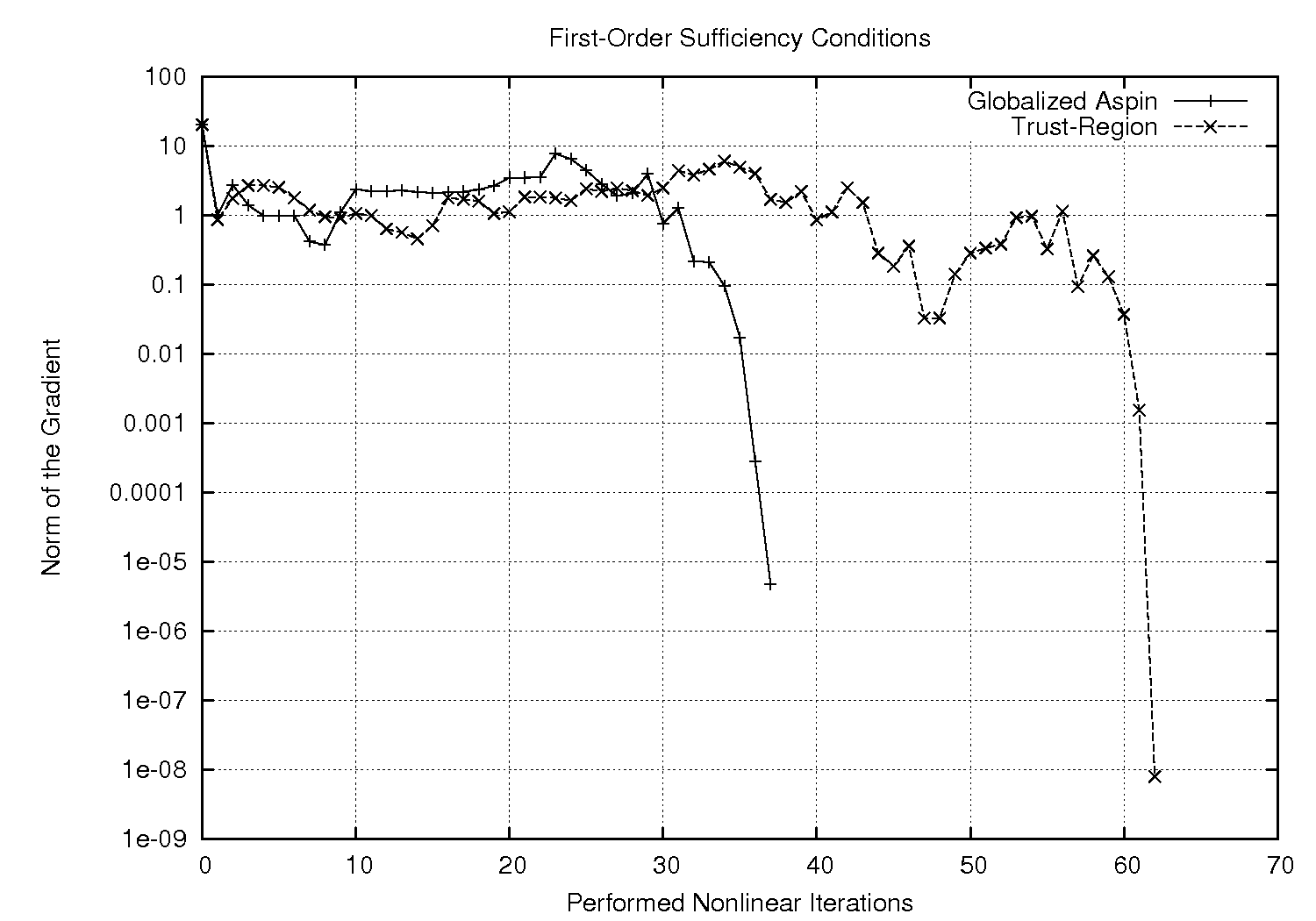}\\
	\ \\ \ \\
	{\bf Computation on 1920 Processors: Trust-Region vs. Globalized ASPIN}\\
	\includegraphics[width=0.45\textwidth]{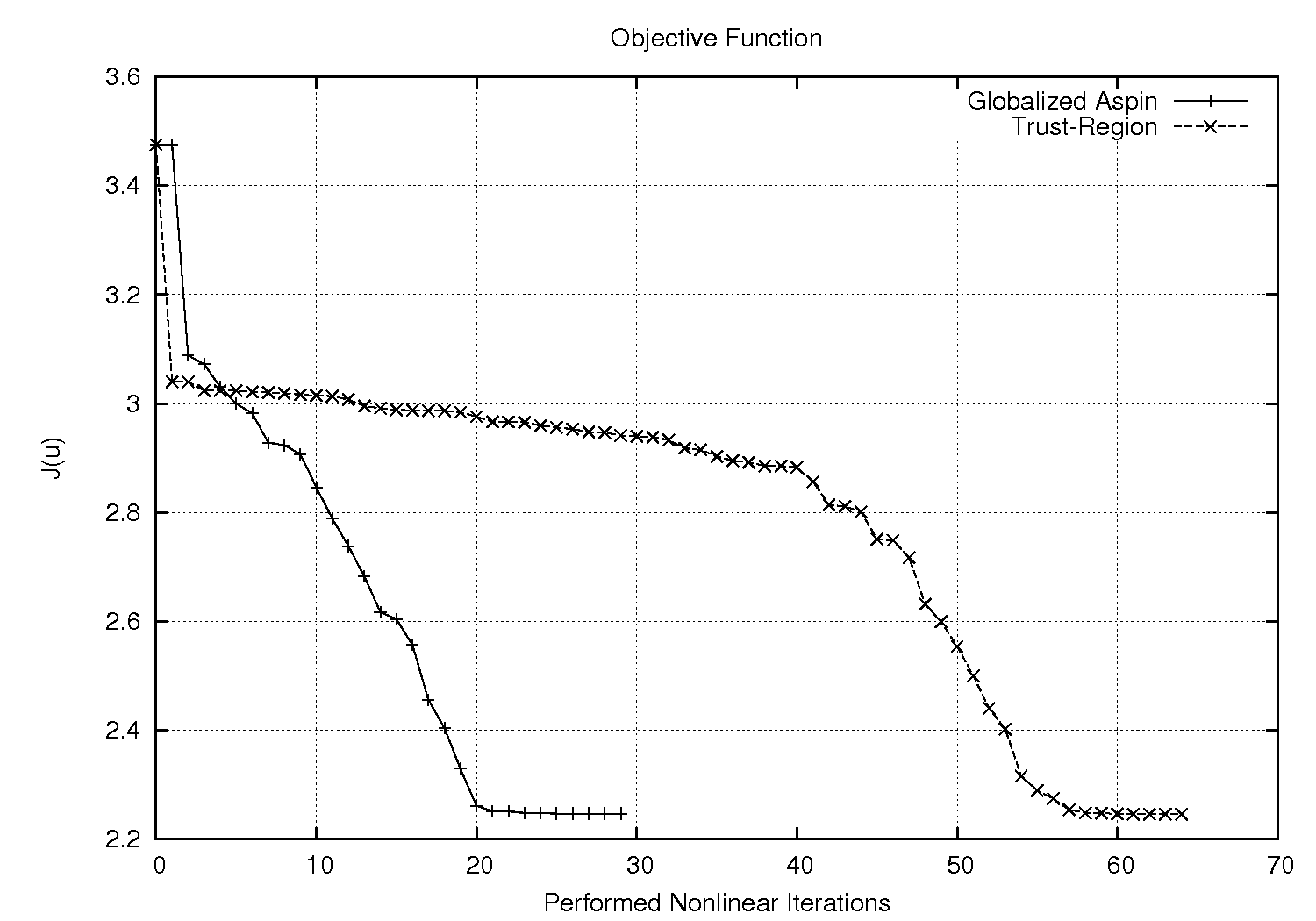} 
	\includegraphics[width=0.45\textwidth]{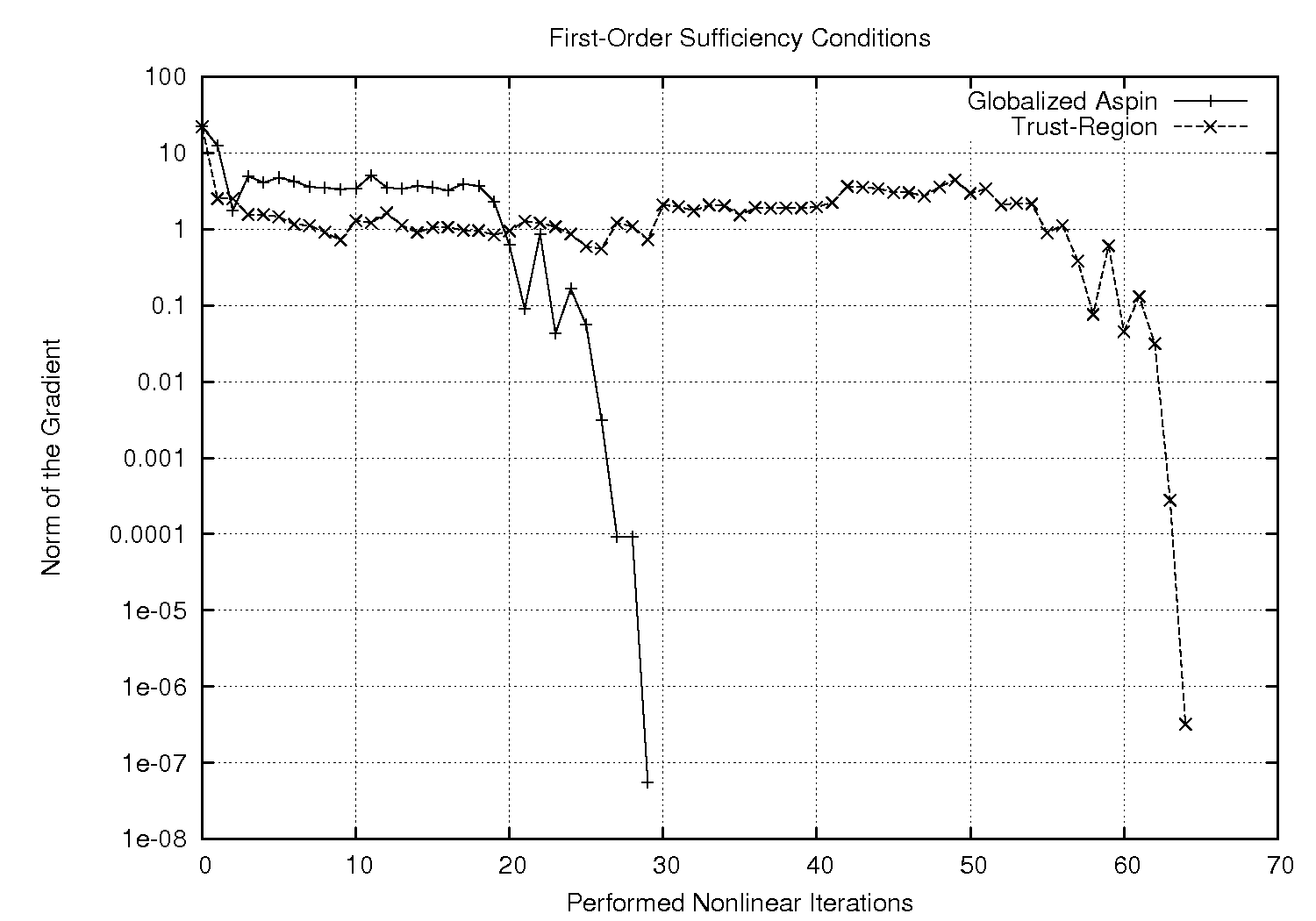}\\
	\ \\ \ \\
	{\bf Comparison Globalized ASPIN with different Processor Numbers}\\
	\includegraphics[width=0.45\textwidth]{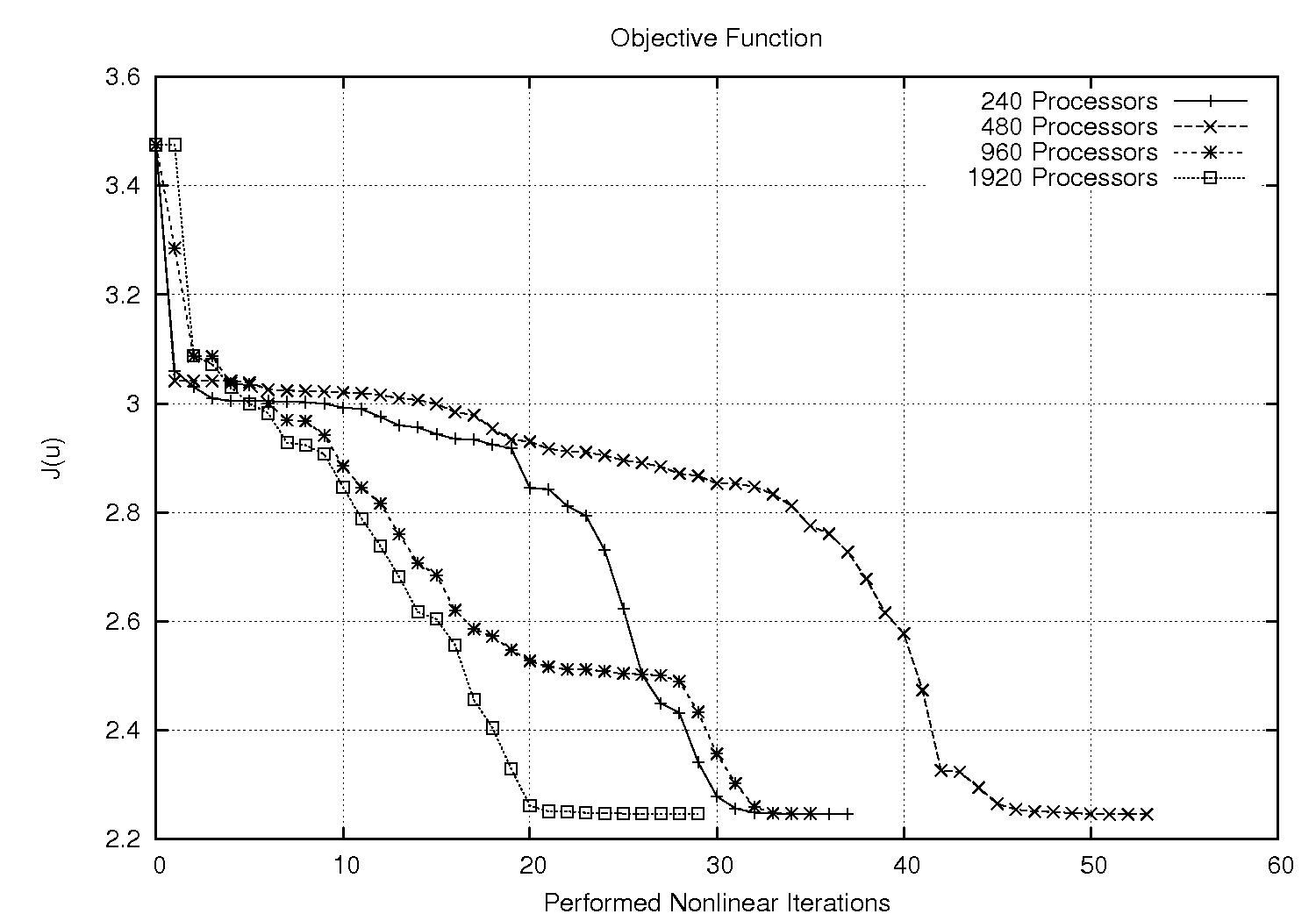}
	\includegraphics[width=0.45\textwidth]{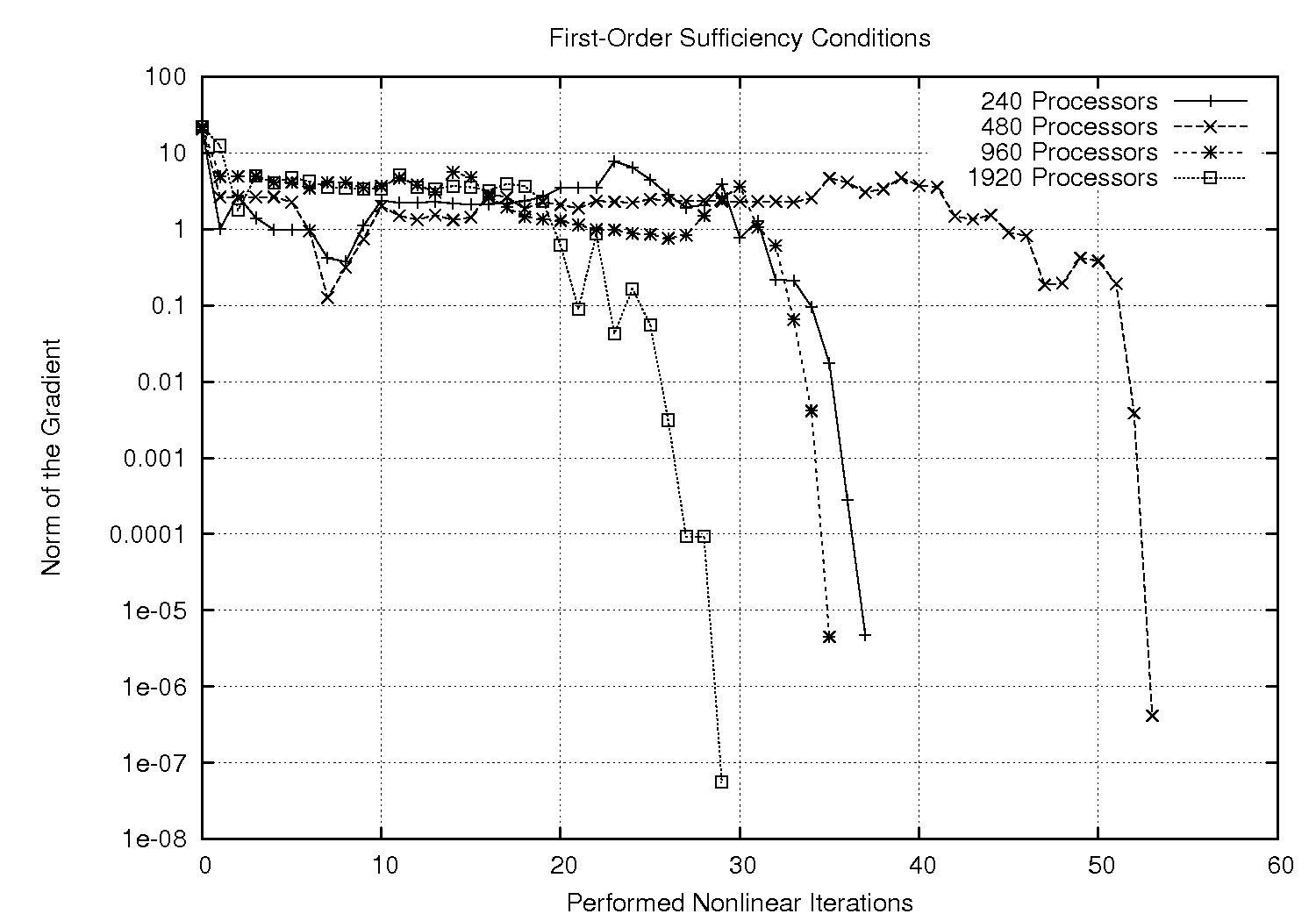}

	 \caption{\label{fig:brickStats}The first-order sufficiency conditions and the values of the objective function 
		after each non-linear iteration for the solution of the example from Section~\ref{sec:brick} with 100,800 unknowns. Computation
		was carried out on $240$, $480$, $960$, $1920$ processors. 
		Here, we compare a trust-region method with the globalized ASPIN approach, where both methods converge to exactly the same
		solution.
	}
\end{figure}		

\begin{table}
	\begin{tabular}{|l|c|c|}
		\hline	
		&Trust-Region & Globalized ASPIN \\		
		\hline
		Overall Time & 3293.95 & 2166.99\\
		Linear Solver for global trust-region problem & 3292.67 & 2042.96 \\
		Linear Solver in local solution phase & --- & 7.61 \\
		Assembling & 12.71 & 68.09 \\
		\hline
	\end{tabular}	
	\caption{\label{table:brickComptimes240} 
		Computation times (in seconds) needed to solve Problem~\eqref{eqn:nlMinProb} on $240$ processors.
		Problem settings are described in Section~\ref{sec:brick}.
	}
\end{table}		

\begin{table}
\begin{tabular}{|l|c|c|}
			\hline	
			&Trust-Region & Globalized ASPIN \\		
			\hline
			Overall Time & 2711.17 & 1238.76\\
			Linear Solver for global trust-region problem & 2722.98 & 1220.27 \\
			Linear Solver in local solution phase & --- & 0.66 \\
			Assembling & 4.91 & 12.40 \\
			\hline
		\end{tabular}
	\caption{\label{table:brickComptimes1920} 
		Computation times (in seconds) needed to solve Problem~\eqref{eqn:nlMinProb} on $1920$ processors.
		Problem settings are described in Section~\ref{sec:brick}.
	}
\end{table}

\begin{table}
	\begin{tabular}{|l|c|c|c|c|}
		\hline	
		&240 cores & 480 cores & 960 cores & 1920 cores\\		
		\hline
		Overall Time & 2166.99 & 2044.14	& 1102.32	& 1238.76\\
		Linear Solver global TR problem & 2042.96 & 1918.69	&	1058.14 &1220.27 \\
		Linear Solver local & 7.61 & 6.65 & 1.77	&0.66 \\
		Assembling & 68.09	 & 71.49	& 23.77 & 12.40 \\
		Nonlinear Iterations & 37 & 53 & 35 & 29\\
		\hline
	\end{tabular}
	\caption{\label{table:brickComparisons} 
	Computation times (in seconds) needed to solve Problem~\eqref{eqn:nlMinProb} with different numbers of processors.
	Problem settings are described in Section~\ref{sec:brick}. The poor scaling from $960$ to $1920$ cores is
	due to the small local problem size. On the other hand, the poor scaling between $240$ and $480$ cores 
	can be explained by an increased number of non-linear iterations.
	}
\end{table}

\subsubsection{Local Solution and Recombination}

For the computation of $\tilde g_i$, we employ the linear recombination approach from Section~\ref{sec:linearApproach}. 
As we pointed out, it is therefore not important which method we use for the local computations on the subsets.
In our implementation, we employ a trust-region approach. 
This trust-region method either stops the local computation if the norm of the local gradient is lower than $10^{-9}$ or after 20 trust-region iterations. 
% In any case, each local trust-region method memorizes the last trust-region radius which allows
% for eventually computing local trust-region corrections.

\subsubsection{Nonlinear Solution Process}

The numerical results of the following sections compare the runtime and convergence of the presented globalized ASPIN strategy
 to a trust-region strategy without the preconditioning step. Both globalization strategies employ a 
parallelized Steihaug Toint conjugate gradient method~\cite{steihaug1983conjugate,toint1981towards}
 for the solution of the arising trust-region subproblems, which employs a multigrid method as
a smoother. In order to be able to prove a sufficient decrease, we only accept the computed
trust-region steps, if their predicted reduction is better than the one of the respective Cauchy step. 

Alternatively, one might employ a Lanczos method (for broad survey see
\cite{ConnGouldToint00}) for the solution of the arising quadratic programming problems.

\subsection{Deformation of a Hollow Brick}
\label{sec:brick}

In this section, we present computational results showing the efficiency of the new approach for non-linear
programming problems in non-linear elasticity. As a reference application, we consider the compression of a brick as
shown in Figure \ref{fig:brick}. The brick itself is a hollow geometry with dimension $(1.0,1.0,2.0)$ and is compressed
by $10.0\%$. The chosen material parameters are $E=3000$Pa and $\nu=0.3$.

The solution of the resulting minimization problem is shown in Figure~\ref{fig:brick}, the values of the objective function 
and the gradients are given in Figure~\ref{fig:brickStats}. Moreover, computation times of the respective strategies are
given in Tables~\ref{table:brickComptimes240}, \ref{table:brickComptimes1920} and \ref{table:brickComparisons}.

Note that in Figure~\ref{fig:brickStats}, the number of non-linear iterations needed to reach a specific reduction in the objective function $J$
sometimes decreases substantially, when the number of processors is increased.
This would not be the expected behavior usual for parallel Schwarz methods for linear problems.
However, we here are not dealing with a parallelization of the Newton method, but
our method is a non-linear, parallel solver, whose behavior changes, when we employ a different domain decomposition, as different and more quadratic models are employed.

\subsection{Deformation of a Hollow Sphere}
\label{sec:sphere}

\begin{figure}
\centering
\includegraphics[width=0.47\textwidth]{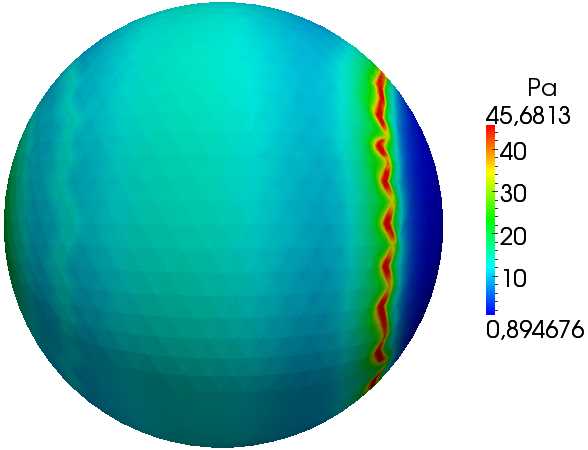}\hspace{0.05\textwidth}
\includegraphics[width=0.35\textwidth]{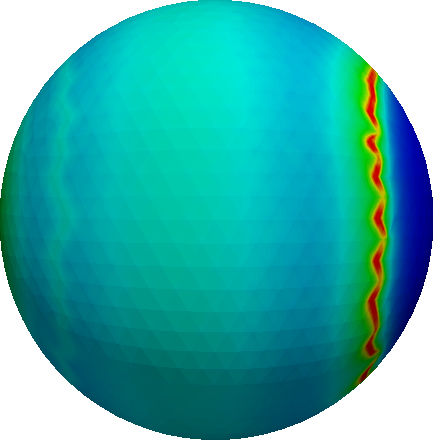}

\caption{\label{fig:sphere}Numerical result of the example from Section \ref{sec:brick} with 828,864 unknowns, where
the initial geometry is deformed according to the solution of \eref{eqn:nlMinProb} (\emph{left figure}). \emph{Right
figure}: the initial geometry, which is -- beside the stresses -- only slightly different from the deformed
configuration.
}
\end{figure}

As second application, we consider the very slight compression of a sphere, as
shown in Figure \ref{fig:sphere}. The sphere is a hollow geometry with inner radius $0.4$ and outer radius $0.5$.
The chosen material parameters are $E=3000$Pa and $\nu=0.3$. 
In this second example, the sphere is fixed at one side, i.e., we apply $0$ displacements at the Dirichlet boundary. On the
opposite side of the geometry, we apply constant forces, i.e., $(40,\ldots,40)^T$ as Neumann values.

\begin{figure}
\centering
{\bf Computation on 32 cores}\\
\includegraphics[width=0.49\textwidth]{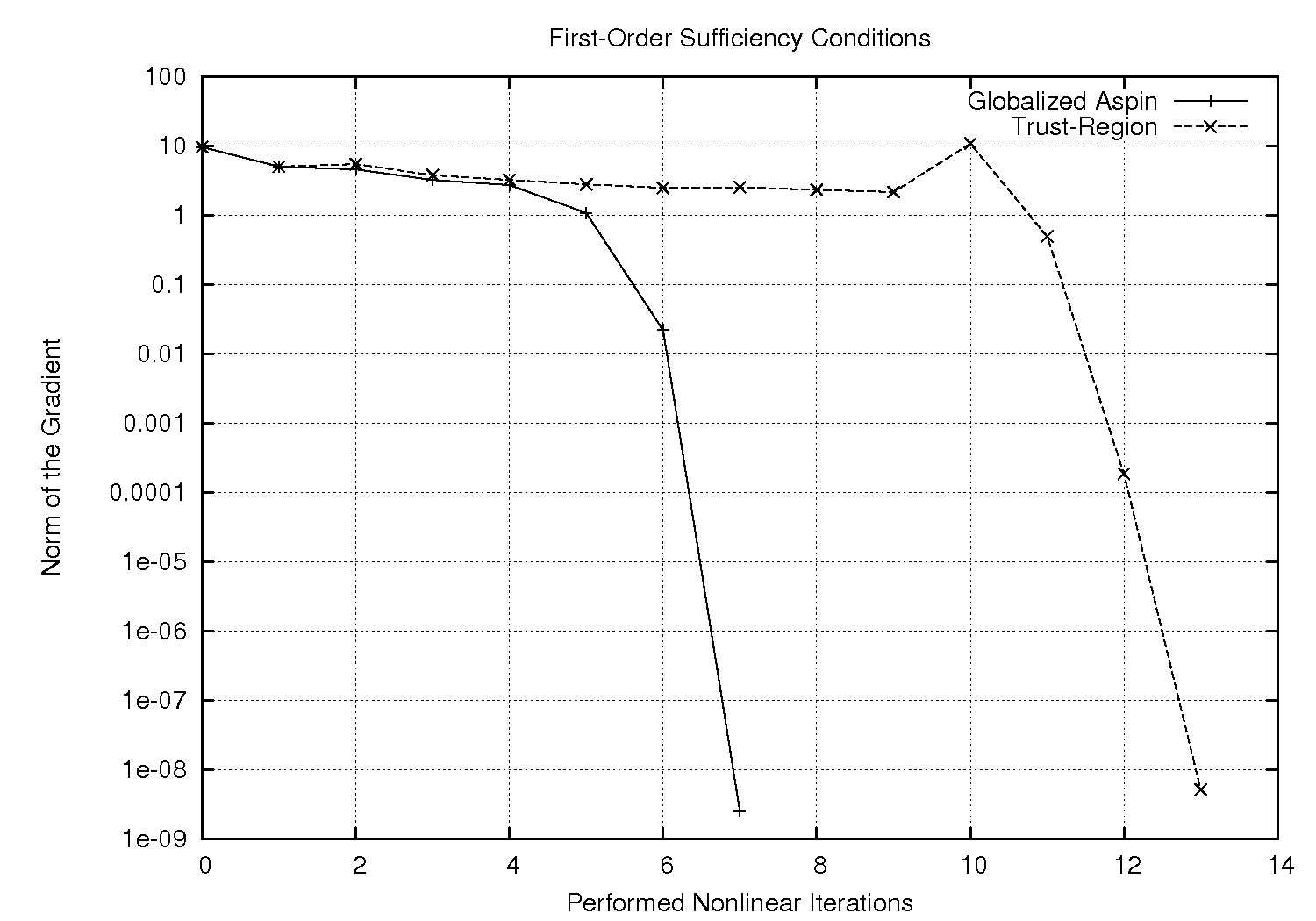}
\includegraphics[width=0.49\textwidth]{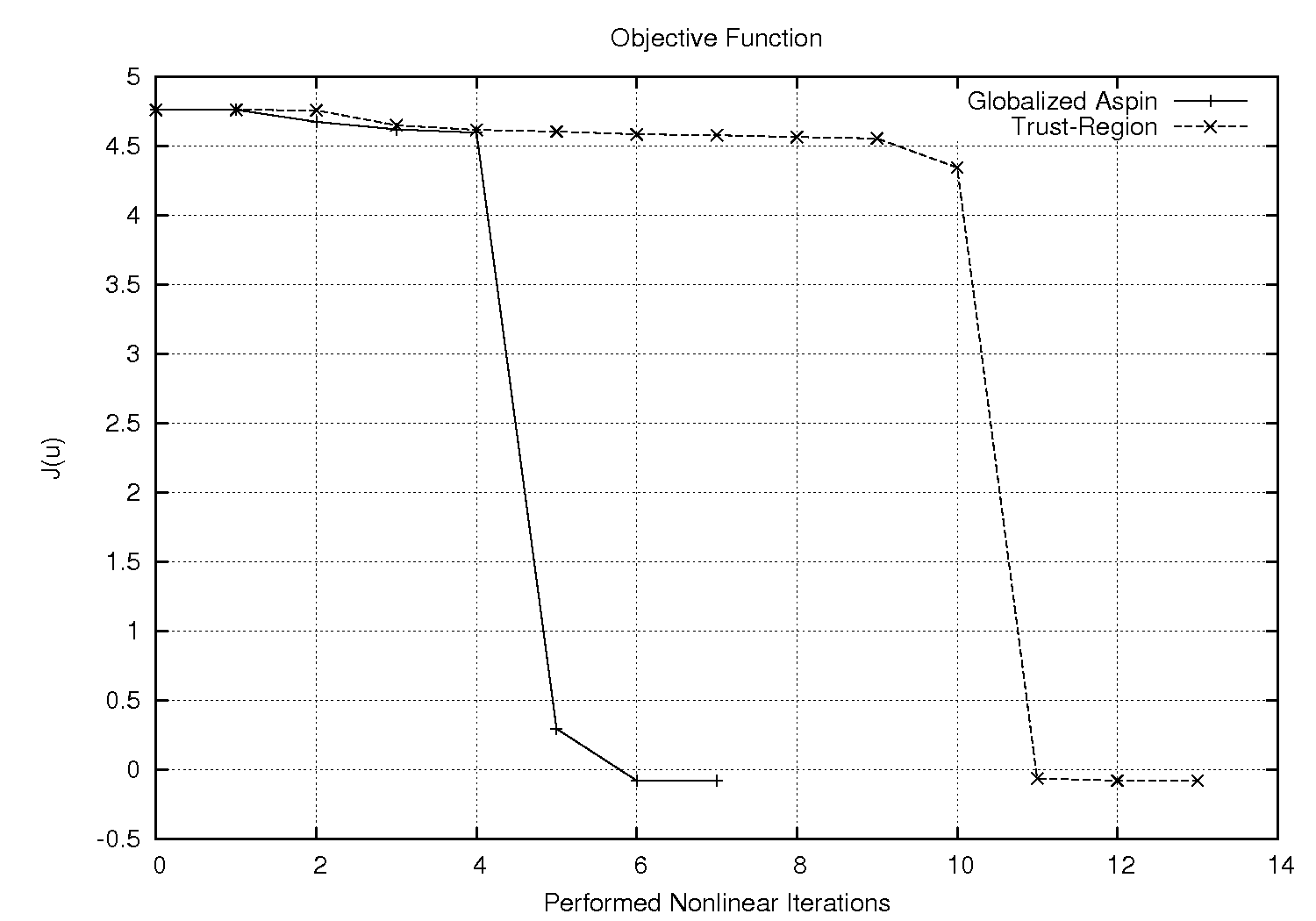}
\
\\[0.3cm]
{\bf Computation on 64 cores}\\
\includegraphics[width=0.49\textwidth]{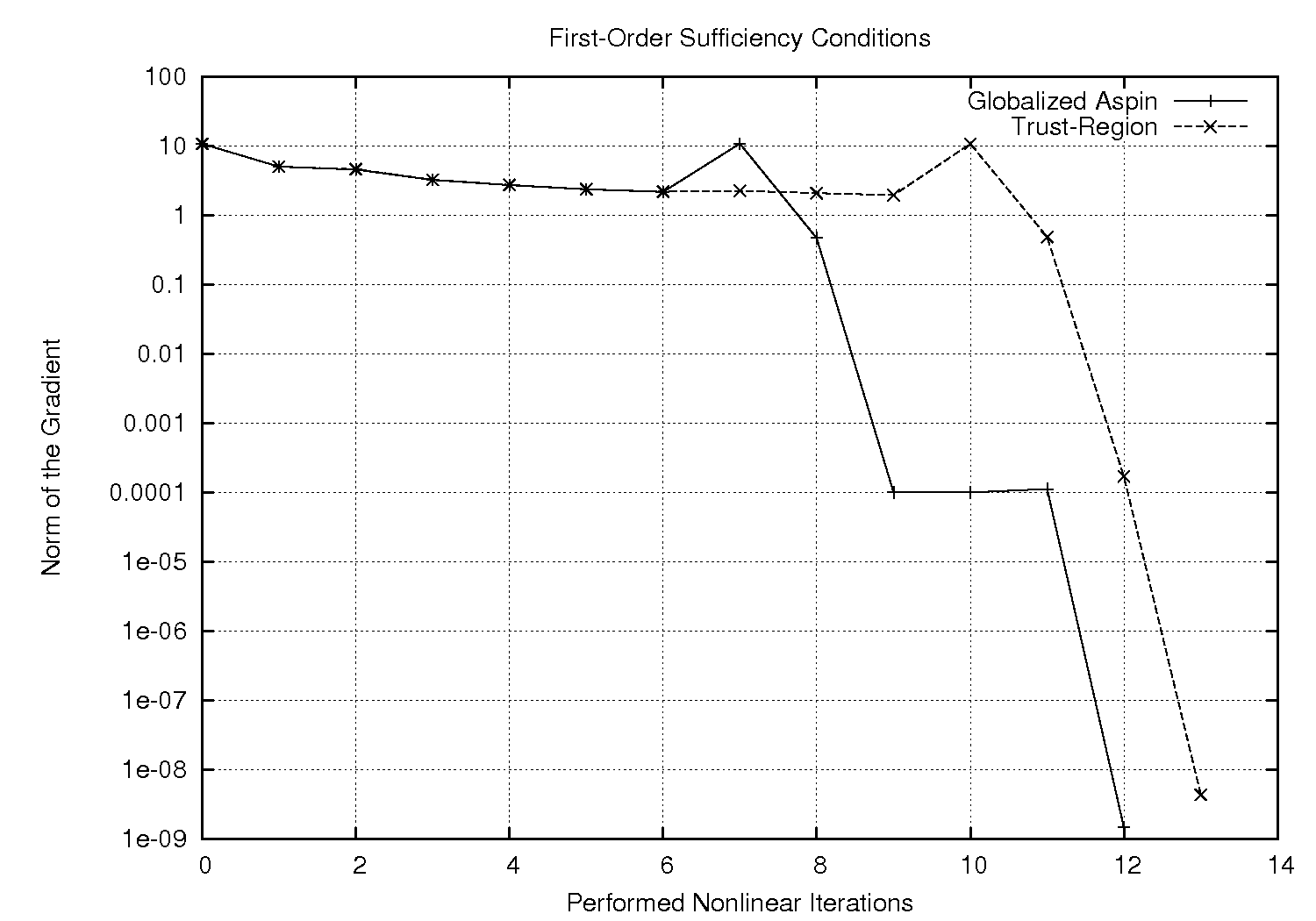}
\includegraphics[width=0.49\textwidth]{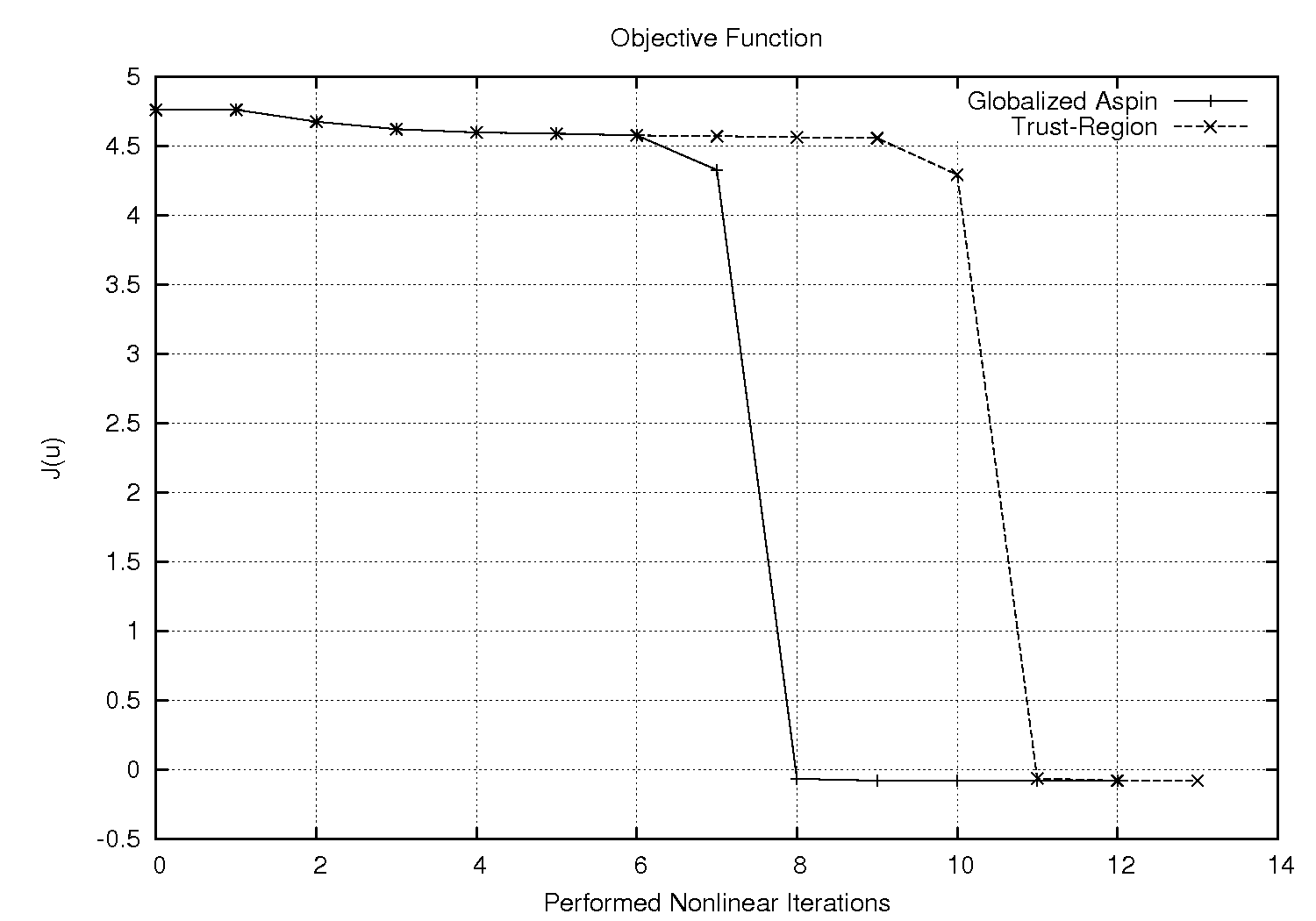}

\caption{\label{fig:sphereStats}The first-order sufficiency conditions  and the values of the objective function
after each non-linear iteration for the solution of the example from Section~\ref{sec:brick} with 828,864 unknowns.
Computation was carried out on 32 and 64 processors. Here, we compare a trust-region method with the globalized ASPIN
approach, where both methods converge to exactly the same solution.
}
\end{figure}

% \subsubsection{Computational Results} 

The solution of the resulting minimization problem is shown in Figure~\ref{fig:sphere}, the values of the objective
function and the gradients are given in Figure~\ref{fig:sphereStats}. Moreover, computation times of the respective
strategies are given in Table~\ref{table:sphereComptimes}. Note that also in this case, all computations end in the
same minimum.

\begin{table}
	\centering
	{\bf Computation on 32 cores}\\
	\begin{tabular}{|l|c|c|}
		\hline	
		& Trust-Region & Globalized ASPIN \\		
		\hline
		Overall Time &3323.96 & 2576.61\\
		Linear Solver for global trust-region problem & 3211.49 & 1796.65 \\
		Linear Solver in local solution phase & --- & 442.23 \\
		Assembling & 70.39 & 299.05 \\
		\hline
	\end{tabular}
	\ \\\ \\\ \\
	{\bf Computation on 64 cores}\\
	\begin{tabular}{|l|c|c|}
		\hline	
		& Trust-Region & Globalized ASPIN \\		
		\hline
		Overall Time &1805.94 & 2211.23\\
		Linear Solver for global trust-region problem & 1748.03 & 1653.55 \\
		Linear Solver in local solution phase & --- & 224.55  \\
		Assembling & 38.80 & 204.65 \\
		\hline
	\end{tabular}
	\caption{\label{table:sphereComptimes} Computation times (in seconds) needed to solve problem
\eref{eqn:nlMinProb} on 64 processors.
	Problem settings are described in Section~\ref{sec:sphere}}.
\end{table}

\section{Conclusions}

In this article, we presented and analyzed a globalization strategy which extends the ASPIN method presented in
\cite{CaiKeyes00}. The key idea of this globalization approach is to consider ASPIN's non-linearly preconditioned Newton
step as the first-order conditions of a particular quadratic programming problem, which we refer to as preconditioned
quadratic model. 

Due to the interpretation of this preconditioned
model as a perturbed model, one can enforce convergence of the method by employing a trust-region algorithm on the one
hand and by controlling the perturbation, on the other. Approaches for controlling perturbations were introduced for
instance in \cite{Toint1988,Carter1993,ConnEtAl1993,Felgenhauer1997}, where the focus is on numerical differentiation
and the solution of constrained non-linear programming problems. Here, however, the perturbation results from
the fact that the gradient is based on the sum of local correction vectors for related, local non-linear programming problems.
A particular feature of the current algorithm is that it exploits (massive) parallelism by solving these local
non-linear programming problems asynchronously and in parallel. In this context, we introduced another trust-region radius, whose purpose is to
control the perturbation of the preconditioned gradient. We introduced two strategies to
compute the preconditioned gradients, which satisfy this perturbation constraint. Furthermore, our numerical results
reveal the stability of the novel method and, in comparison to the trust-region method, a significant speed-up.

\section*{Acknowledgements}
The authors thank Lea Conen for her detailed corrections, 
in particular for rewriting and clarifying the reasoning in  
sections \ref{sec:aspin}, \ref{sec:trustRegion}, \ref{sec:localNonlinearSolution}, and \ref{sec:firstOrderConvergence}.

\bibliographystyle{plain}
\bibliography{gaspin}

% \clearpage
% \input{comments.tex}

\end{document}

%% file: commands.tex
\providecommand{\abs}[1]{\left\lvert#1\right\rvert}
\providecommand{\norm}[1]{\left \lVert#1\right\rVert}

\newcommand{\globiti}{\ensuremath{u_{i}}}
\newcommand{\globitstart}{\ensuremath{u_{0}}}
\newcommand{\globstepi}{\ensuremath{s_{i}}}

\newcommand{\localu}{\ensuremath{u^k}}
\newcommand{\locals}{\ensuremath{s^k}}

\newcommand{\Pk}{\ensuremath{P^{k}}}
\newcommand{\Rk}{\ensuremath{R^k}}
\newcommand{\Ik}{\ensuremath{I^k}}
\newcommand{\Xk}{\ensuremath{X^k}}
\newcommand{\Hk}{\ensuremath{H^k}}
\newcommand{\Jk}{\ensuremath{J^k}}
\newcommand{\Bk}{\ensuremath{B^k}}
\newcommand{\Ck}{\ensuremath{C^k}}

\newcommand{\Sk}{\ensuremath{\D^k}}

\DeclareMathOperator*{\argmin}{argmin}

\newcommand{\ignore}[1]{}
\newcommand{\nobibentry}[1]{{\let\nocite\ignore\bibentry{#1}}}

%% file: gaspin_definition.tex
\section{A Globalized ASPIN Framework}
\label{sec:localNonlinearSolution}

If the global Newton step is solved in direction of the preconditioned ASPIN gradient,
the resulting step does not necessarily lead to a sufficient decrease in the objective function,
and the method does not necessarily converge for starting vectors far away from first-order critical points.  
In the next step, we thus combine the trust-region method with the ASPIN ideas
in order to define a global, parallel solution method for non-linear optimization problems.
That is, we want in particular to answer the question 
how to compute global corrections~$\globstepi$ based on the local corrections~$\locals$ from \eref{eqn:aspinStep} such that the resulting algorithm is a globalization strategy.

For that purpose, we modify the trust-region algorithm, Algorithm~\ref{alg:TrustRegionSolver}.
In particular, we make the following changes:
\begin{itemize}
	\item In Problem~\eqref{eqn:localMinProb}, we exchange the quadratic model function $\psi_i$ by a modified version~$\tilde \psi_i$, the \emph{preconditioned quadratic model}. 
		It additionally depends on the product of the modified ASPIN gradient and the inverse of the non-linear Schwarz preconditioner $C_i\cdot g_i^{\text{ASPIN}}$. 
% 		That is, we will interpret the preconditioned Newton system in  \eref{eqn:aspinNewtonStep}
% 		as the solution of a perturbed quadratic programming problem~\eqref{eqn:globalQP}, compare Remark~\ref{rem:relationMinProbSystem}.
		We have to introduce strategies in order to control this perturbation.
% 		These modifications necessary for global convergence then lead to the function~$\tilde \psi_i$.
	\item As common in trust-region strategies, the modified quadratic function~$\tilde \psi_i$ will also be used to approximate the actual reduction. 
		Therefore, the definition of the decrease ratio~$\rho_i$ is adapted.
	\item Furthermore, we modify the sufficient decrease condition in \eref{eqn:sufficientDecrease}, in order to incorporate the additional quantities. 
\end{itemize}
The following sections then deal with how to exactly make these changes such that the resulting algorithm converges for arbitrary starting vectors 
and employs ASPIN inspired strategies.

\subsection{The Preconditioned Quadratic Model}
\label{sec:preconditionedQuadraticModel}
In this section, we explain how we compute global corrections~$\globstepi$ based on the local ASPIN steps~$\locals$.
For that purpose, we modify the quadratic model~$\psi_i$ from \eref{eqn:localQuadraticModel}
in order to include the local ASPIN corrections defined in \eref{eqn:aspinStep} 
and the additive Schwarz preconditioner~$C_i^{-1}$ defined in \eref{eqn:matrixC}.
That is, we introduce the following, \emph{preconditioned quadratic programming problem}
\begin{equation}
\label{eqn:globalQP}
	\globstepi := \argmin_{s \in \mathbb{R}^n,\,\norm{s} \le \Delta_{i}^G} \widetilde \psi_{i}(s)
\end{equation}
where we employ the following \emph{preconditioned quadratic model}
\begin{equation}
\label{eqn:globalTrModel}
	\widetilde \psi_i(s) = \frac{1}{2}\langle s,   B_i s \rangle + \langle \tilde g_i, s \rangle. 
\end{equation}
The matrix $B_i=B(\globiti)$ is, as before, a symmetric approximation to the Hessian $\nabla^2 J\left(u_i\right)$.
The vector $\tilde g_i \in \mathbb{R}^n$ will be referred to as the \emph{preconditioned gradient}. 
Its definition is based on the original gradient~$g_i = \nabla J\left(\globiti\right)$ and on the product of ASPIN related quantities~$C_i\cdot g_i^{\text{ASPIN}}$.
That is, it can be defined as
\begin{equation*}
	\tilde g_i = h(g_i, C_i \cdot g_i^{\text{ASPIN}}),
\end{equation*}
where the function $h: \mathbb R^n \times \mathbb R^n \to \mathbb R^n$
combines the original gradient and ASPIN's gradient in a suitable way, 
such that the conditions specified in Section~\ref{sec:assumptionsOnGradient} are fulfilled.
These conditions are needed in order to ensure that 
the modified model~$\widetilde \psi_i$ is a reasonably good approximation to the actual reduction~$\ared_i$.
Two different approaches for actually implementing $h$, i.e., for computing $\tilde g_i$, will be presented in 
Section~\ref{sec:trustRegionApproach} and Section~\ref{sec:linearApproach}.

We make two remarks in order to ease the understanding of the definition of the modified model and its relation to the original ASPIN method.
First, the modified quadratic model~$\widetilde \psi_i$ gives rise to an ASPIN step under certain conditions.
In fact, by the result cited in Remark~\ref{rem:relationMinProbSystem}, 
if the minimizer~$\globstepi$ of \eref{eqn:globalQP} lies in the interior of the trust region, 
that is, if $\norm{\globstepi} < \Delta_i^G$,
and if $B_i$ is a positive semi-definite matrix,
$s_i$ is the solution of the linear systems of equations
\begin{equation*}
	B_i s_i = -\tilde g_i.
\end{equation*}
If furthermore $\tilde g_i = C_i \cdot g_i^{\text{ASPIN}}$ and assuming that $C_i^{-1}$ exists, 
this system can be reformulated as 
\begin{equation*}
	C_i^{-1} B_i s = - g_i^\text{ASPIN},
\end{equation*}
which is the original ASPIN step~\eqref{eqn:aspinNewtonStep} using an approximated Hessian matrix $B_i \approx \nabla^2 J\left(\globiti\right)$.

Second, we remark that we multiplied $C_i$ with $g_i^\text{ASPIN}$ in \eref{eqn:globalTrModel} instead of using $C_i^{-1}$ separately as in the original ASPIN method. 
The reason for this is that in this way we are able to control the perturbation by just imposing a condition on the gradient,
as we have included all the ASPIN related quantities in the modified gradient in \eref{eqn:globalTrModel}.
Note however that the computation of $C_i$ might be expensive, 
depending on the chosen domain decomposition and on the matrices $\Bk(\Pk\globiti)$.
In the following cases, the vector $C_i \cdot g_i^{\text{ASPIN}}$
can efficiently be computed:
\begin{itemize}
	\item If a non-overlapping domain decomposition is employed,  
	we obtain a non-overlapping block structure of
	$\sum_k \Ck$ which gives rise to
	\begin{equation*}
	C_i = \sum_k \Ck = \sum_k \left(\Ik \Bk(\Pk \globiti) \Rk\right).
	\end{equation*}
	Note that by ``non-overlapping'' we here mean that no two subdomains share any degrees of freedom, not even on the interface.
\item As long as $(C_i)^{-1}$ is a sparse matrix, one might employ
	a preconditioned Krylov method in order to compute~$C_i g_i^\text{ASPIN}$.
\end{itemize}

\subsection{Assumptions on the Preconditioned Gradient}
\label{sec:assumptionsOnGradient}

% For the analysis in Section~\ref{sec:firstOrderConvergence}, 
% it is not necessary to state 
% \todolea{state?}
% further assumptions directly on the local solution process. 
% Therefore, in this section we will state assumptions on the preconditioned gradient, since 
% computing global Newton steps (by solving \eref{eqn:globalQP}) 
% and summing them up does not suffice for proving convergence of the resulting algorithm. 

In order for the trust-region method employing the preconditioned quadratic model~$\widetilde \psi_i$ to converge -- or to be able to prove convergence -- 
it is necessary to control the behavior of the preconditioned gradient~$\tilde g_i$.
In particular, we assume that the preconditioned gradient satisfies
\begin{equation}
\label{eqn:assumptionDeltaS}
	\|\tilde g_i - g_i\| \le \Delta_i^L.
\end{equation}
Here $\Delta_i^L >0$ is a second trust-region radius used to control the (local) perturbation in the gradient. 
Note that this is in contrast to \cite{Toint1988,ConnEtAl1993}, where just one trust-region radius is employed.
% Condition~\eqref{eqn:assumptionDeltaS} implies that
% \begin{equation}
% \label{eqn:assumptionDeltaS2}
% 	\|g_i\| -\Delta_i^L \le \|\tilde g_i\| \le \|g_i\| + \Delta_i^L. 
% \end{equation}

Condition \eqref{eqn:assumptionDeltaS} plays a crucial role for the analysis of the 
non-linearly left preconditioned trust-region algorithm in Section~\ref{sec:firstOrderConvergence}. 
In fact, this property will enable us to control the perturbation in case that the preconditioned 
gradient does not yield a sufficient decrease. In particular, we have the following 
first-order consistency relationship
\begin{equation*}
	\Delta^L_i \to 0 \qquad \Rightarrow \qquad \tilde g_i \to g_i.
\end{equation*}
In Sections~\ref{sec:trustRegionApproach} and \ref{sec:linearApproach}
we will consider simple and implementable strategies to compute a perturbed gradient $\tilde g_i$ which 
satisfies \eref{eqn:assumptionDeltaS}, while using the local ASPIN corrections whenever suitable.

\subsection{Trust-Region Update and Sufficient Decrease Conditions}
\label{sec:trUpdateSuffDec}

In order to show convergence, we will assume that an extended sufficient decrease condition of the following kind holds:
\begin{subequations}
\label{eqn:modifiedSuffDecrease}
\begin{align}
\label{eqn:modifiedSuffDecrease1}
	- \widetilde \psi_i(\globstepi) 
		&\ge 
		\beta_1 \min\left\{\|\tilde g_i\|^2, \|\tilde g_i\| \Delta_i^G \right\} \\
\label{eqn:modifiedSuffDecrease2}
		&\ge
		\beta_2 \min\left\{\|g_i\|^2, \|g_i\| \Delta_i^G \right\},
\end{align}
\end{subequations}
where $\beta_1 > \beta_2>0$. The first inequality \eqref{eqn:modifiedSuffDecrease1} 
can in general be satisfied by computing 
the Cauchy point, cf.\ the remark after \eref{eqn:sufficientDecrease}.
On the other hand, for arbitrary local corrections~$\locals$, Condition~\eqref{eqn:modifiedSuffDecrease2} generally does not hold. 
We will comment on the actual computation of this extended Cauchy condition in Section~\ref{sec:numerics}.
Therefore, as we have seen, we will handle two trust-region radii, 
$\Delta_i^G$ and $\Delta^L_i$, where the update for $\Delta^L_i$ is based on \eqref{eqn:modifiedSuffDecrease2},
that is
\begin{equation*}
	\beta_1 \min\left\{\|\tilde g_i\|^2, \|\tilde g_i\| \Delta_i^G \right\} \\
	\ge
	\beta_2 \min\left\{\|g_i\|^2, \|g_i\| \Delta_i^G \right\}.
\end{equation*}
In particular, we define an \emph{intermediate} radius as 
\begin{equation}
\label{eqn:deltaLUpdate}
	\tilde \Delta^{L}_{i+1} = 
	\begin{cases}
		\gamma_2 \Delta^L_i & \text{ if \eqref{eqn:modifiedSuffDecrease2} holds }\\
		\gamma_1 \Delta^L_i & \text{ otherwise,}
	\end{cases}
\end{equation}
where $\gamma_1 \in (0,1)$ and $\gamma_2 > 1$. 
Having the intermediate radius $\tilde \Delta_{i+1}^L$, the new local radius~$\Delta^{L}_{i+1}$ is given as 
\begin{equation}
\label{eqn:boundDeltaL}
	\Delta^{L}_{i+1} = \min\left\{\Delta^G_i, \tilde \Delta^L_{i+1}\right\}.
\end{equation}
On the other hand, $\Delta_i^G$ is updated employing the following equation
\begin{equation}
\label{eqn:deltaGUpdate}
	\Delta^{G}_{i+1} = 
	\begin{cases}
		\gamma_2 \Delta^G_i & \text{ if } \tilde \rho_i \ge \eta \text{ and \eqref{eqn:modifiedSuffDecrease2} hold} \\
		\Delta^G_i & \text{ if } \tilde \rho_i \ge \eta \text{ holds} \\
		\gamma_1 \Delta^G_i & \text{ otherwise,}
	\end{cases}
\end{equation}
where
\begin{equation}
\label{eqn:decreaseratio}
	\tilde \rho_i = \frac{\ared_i(\globstepi)}{-\widetilde \psi_i(\globstepi)}
\end{equation}
is the modified version of the decrease ratio in \eref{eqn:decreaseRatio} and $\eta\in(0,1)$.
The global corrections will only be applied if both the sufficient decrease condition in Equation~\eqref{eqn:modifiedSuffDecrease2} and $\tilde \rho_i \ge \eta$ hold, i.e.,
\begin{equation}
\label{eqn:updateU}
	u_{i+1} = 
	\begin{cases}
		\globiti + \globstepi & \text{ if } \tilde \rho_i \ge \eta \text{ and \eqref{eqn:modifiedSuffDecrease2} hold,} \\
		\globiti & \text{ otherwise.}
	\end{cases}
\end{equation}
A feature of the presented algorithm is that we stall the solution process in the global context, 
as long as $\tilde g_i$ does not satisfy Condition~\eqref{eqn:modifiedSuffDecrease2}. 
This approach is reasonable, since we want to distinguish between two different error sources: 
\begin{itemize}
	\item the approximation strength of the trust-region model as a Taylor approximation to the actual decrease, and
	\item the perturbation of the employed preconditioned gradient $\tilde g_i$ in the preconditioned model.
\end{itemize}
This means that if it is certain that the local solution process does not yield gradients $\tilde g_i$
which satisfy Inequality~\eqref{eqn:modifiedSuffDecrease2}, the $\Delta_i^L$ will be reduced which yields $\tilde g_i\to g_i$. 
If the Taylor approximation is poor, we reduce both $\Delta_i^G$ and $\Delta_i^L$ in order to increase $\tilde \rho_i$. 
Summing up these steps yields Algorithm~\ref{alg:gaspin}.

\begin{algorithm}[t]
	\caption{Globalized ASPIN Strategy}
		\label{alg:gaspin}
		\begin{algorithmic}[1]
		\Statex
		\Input{
			$\globitstart \in \mathbb{R}^n$,
			$J:\mathbb R^n \to \mathbb R$, $\D_1,\ldots,\D_N\subset \mathbb R^n$, $\Delta^L_0$, $\Delta^G_0$,
			$\gamma_1, \eta\in (0,1)$, $\gamma_2>1$
		}
		\Output{
			Final iterate $\globiti$
		}
		\Statex
		\State $i \gets 0$
			\While {not converged or maximum number of iterations not reached}
				\State 
				\parbox[t]{\dimexpr\linewidth-\algorithmicindent}{
					On each subdomain $\Sk$ compute a local ASPIN correction~$\locals$ as defined in \eref{eqn:aspinStep}.
					\strut
				}
				\State 
				\parbox[t]{\dimexpr\linewidth-\algorithmicindent}{
    				Compute $\tilde g_i$ based on the subspace corrections~$\locals$,
					such that it satisfies \eref{eqn:assumptionDeltaS} (cf., Sections~\ref{sec:linearApproach} and~\ref{sec:trustRegionApproach}).
    			\strut} 
				\State Solve \eref{eqn:globalQP} in order to obtain a global correction $\globstepi$ satisfying \eqref{eqn:modifiedSuffDecrease1}.
				\State Compute $\Delta^L_{i+1}$ and $\Delta^G_{i+1}$ according to equations~\eqref{eqn:deltaLUpdate} and \eqref{eqn:deltaGUpdate}.
				\State Update $u_{i+1}$ by means of \eref{eqn:updateU}. 
				\State $i\gets i+1$
			\EndWhile
			\State \Return $\globiti$
	\end{algorithmic}
\end{algorithm}

\begin{remark}
\label{rem:perturbedDecreaseRatio}
As a crucial ingredient for the analysis of trust-region methods one exploits that for successful steps
\begin{equation*}
	\ared(\globstepi) \ge \eta\;\pred(\globstepi) \ge \beta_2 \min\left\{\|g_i\|^2, \|g_i\| \Delta_i^G \right\}
\end{equation*}
holds.
In our context, we do not employ the quadratic programming problem in \eref{eqn:localMinProb} for computing 
trust-region corrections, but the perturbed one in \eref{eqn:globalQP}. 
This gives rise to a perturbed result for successful steps, i.e.,
\begin{equation*}
	\ared(\globstepi) \ge -\eta \tilde \psi_i(\globstepi) \ge \beta_1 \min\left\{\|\tilde g_i\|^2, \|\tilde g_i\| \Delta_i^G \right\}.
\end{equation*}
As pointed out before, the Cauchy point satisfies this condition. Then,
the conditions \eqref{eqn:modifiedSuffDecrease1} and \eqref{eqn:modifiedSuffDecrease2} give for successful steps 
\begin{equation*}
	\ared(\globstepi) \ge -\eta \tilde \psi_i(\globstepi) \ge \beta_2 \min\left\{\|g_i\|^2, \|g_i\| \Delta_i^G \right\}.
\end{equation*}
% Therefore, condition \eref{eqn:modifiedSuffDecrease} allows for controlling the perturbation of the 
% trust-region model equation \eref{eqn:globalQP} and, in return, for proving global convergence 
% of the left preconditioned trust-region strategy as in the analysis in Section~\ref{sec:firstOrderConvergence}.
\end{remark}

\subsection{Particular Strategies for the Computation of the Preconditioned Gradient}
\label{sec:particularStrategies}

In Section~\ref{sec:assumptionsOnGradient}, we have stated particular assumptions on the preconditioned gradient.
Now, we will introduce two approaches, a trust-region and a damping approach, in order to combine the asynchronously computed
subspace corrections to a preconditioned gradient which satisfies \eref{eqn:assumptionDeltaS}. Prior to the presentation 
of the approaches in Sections~\ref{sec:trustRegionApproach} and \ref{sec:linearApproach}, we remark that 
if $\Delta_i^L$ and $\Delta_i^G$ are sufficiently large,
if all corrections computed in \eref{eqn:globalQP} are accepted and lie in the interior of the trust region,
and if the approximation to the Hessian is positive semidefinite, 
our globalized ASPIN strategy reduces to the ASPIN method 
presented in Section~\ref{sec:aspin}, cf.\ the remark in Section~\ref{sec:preconditionedQuadraticModel}.

\subsubsection{A Trust-Region Approach}
\label{sec:trustRegionApproach}

Here, we derive an approach to directly control the perturbation of the gradient on the 
subsets by stating a constraint on the length of the local corrections~$\locals$.
% This approach stands in contrast to Section~\ref{sec:linearApproach}, where we employ a posteriori control
% the length of the local corrections. 
For this approach, the local solution process,
that is, the computation of the vectors~$\locals$ in \eref{eqn:aspinStep},
 is assumed to start from 
\begin{equation}
\label{eqn:localStart}
	\Pk \globiti - \left(\Bk(\Pk \globiti)\right)^{-1} \nabla \Hk(\Pk \globiti) \in \Sk,
\end{equation}
where the projected global iterate~$\Pk \globiti$ serves, as before, as the initial iterate on the subset.
Therefore, after projecting the current global iterate to the subset, the first Newton step will be computed and fully accepted
by definition of the initial iterate.

We note that for the definition of the perturbed model~$\widetilde \psi_i$, 
we used $C_i$ instead of $C_i^{-1}$ and thus -- in case of a non-overlapping decomposition -- only had to assume
that $\Bk(\Pk \globiti)$ exists, which is weaker than assuming that $\Bk(\Pk \globiti)^{-1}$ exists. 
However, for the computation of \eref{eqn:localStart} in the trust-region approach presented in this section, we have to assume that also the inverse exists.
This is for instance the case if $\Bk$ is defined by the BFGS method, cf.~\cite{NocedalWright2006}.

We assume that the local computation, which is the application of a limited number of
trust-region steps, then produces the subset correction
$$\locals = \localu_f-\Pk \globiti,$$ 
where $\localu_f\in \Sk$ denotes the final iterate of the subset computation.

Furthermore, in order to satisfy \eref{eqn:assumptionDeltaS}, we assume that 
\begin{equation}
\label{eqn:localConstraint}
	\norm{\locals + \left(\Bk(\Pk \globiti)\right)^{-1} \nabla \Hk(\Pk \globiti)} \le \omega \Delta_i^L
\end{equation}
holds, where 
$$
	0<\omega\le\frac{N}{\|C_i\|}\Delta_i^L.
$$
This means, that the local trust-region steps~$\locals$ are not allowed to move further
away from the initial Newton step than $\omega \Delta_i^L$.

In order to compute a correction which satisfies \eref{eqn:localConstraint}, a modified version
of Algorithm~\ref{alg:TrustRegionSolver} can be employed. The necessary change is small.
Here, we follow the local trust-region approach 
in \cite{GrattonSartenaerToint06} and modify the trust-region update as follows. 
In the $l$-th iteration, 
we compute the $l$-th local trust-region step~$\locals_l$ on $\Sk$.
We employ \eref{eqn:updateTrustRegionRadius} to compute an intermediate \emph{local} radius~$\widetilde \Delta^k_{l+1}$ for the next iteration on $\Sk$. 
Then, we choose the actual trust-region radius for the next trust-region step as 
\begin{equation*}
	\Delta^k_{l+1} = \min\left\{\widetilde \Delta^k_{l+1}, \Delta^L_i - \|\locals_{l} - \left(\Bk(\Pk \globiti)\right)^{-1} \nabla
\Hk(\Pk \globiti)\|\right\}.
\end{equation*}

As shown, e.g., in \cite[Lemma 2.1]{GrattonSartenaerToint06}, this local trust-region algorithm
computes a local correction $\locals = - \left(\Bk(\Pk \globiti)\right)^{-1} \nabla \Hk(\Pk \globiti) + \sum_{l=1}^{n_L} \locals_{l}$ 	
satisfying \eref{eqn:localConstraint}. 
Here, $n_L\in \mathbb N$ denotes the number of local trust-region steps.

Now, we define the preconditioned gradient in \eref{eqn:globalTrModel} as
\begin{equation}
\label{eqn:trUpdateStep}
	\tilde g_i = -C_i \cdot \sum_k \Ik \locals,
\end{equation}
where each subset correction~$\locals$ is given by
\begin{equation}
\label{eqn:localCorrection}
	\locals = \localu_{f} - \Pk\globiti = -\left(\Bk(\Pk \globiti)\right)^{-1} \nabla \Hk(\Pk \globiti)+ \sum_{i=1}^{n_L} \locals_{i}.
\end{equation}

As pointed out, this trust-region approach under certain conditions gives rise to an ASPIN Newton step, cf.\ \eref{eqn:aspinNewtonStep}. 
Furthermore, the following lemma shows that the additional trust-region constraint yields that Assumption~\eqref{eqn:assumptionDeltaS} will be satisfied by the computed corrections.

\begin{lem}
\label{lem:sufficientDecreaseTrustRegionApproach}
	Let assumptions (\ATR1), (\ATR2) and (\ATR3) hold.
	Then the preconditioned gradient, computed as the 
	trust-region update step \eqref{eqn:trUpdateStep}, satisfies \eref{eqn:assumptionDeltaS}.	
\end{lem}
\begin{proof}	
	We consider the preconditioned gradient $\tilde g_i$ given by \eref{eqn:trUpdateStep}.
	Due to the definition of the local correction in \eref{eqn:localCorrection}, 
	it is natural to split it as follows
	\begin{equation*}
		\locals = -\left(\Bk(\Pk \globiti)\right)^{-1} \nabla \Hk(\Pk \globiti) + \delta \locals,
	\end{equation*}
	where we abbreviate $\delta \locals = \sum_{i=1}^{n_L} \locals_{i}$.
	
	Furthermore, by \eref{eqn:exactLowerLevelModel2}
	we have that $\nabla \Hk(\Pk \globiti) = \Rk g_i$.
	This gives rise to
	\begin{subequations}
	\begin{align}
		\label{eqn:representation1}
		\tilde g_i &= C_i\cdot\sum_k\left( \left[\Ik \left(\Bk \right)^{-1} \Rk\right] g_i - \Ik \delta \locals \right)\\
			&= C_i\cdot (C_i)^{-1} g_i - C_i\sum_k \Ik \delta \locals )\\ 
			&=  g_i - C_i\cdot\sum_k  \Ik \delta \locals.
	\end{align}
	\end{subequations}		
	Exploiting \eref{eqn:localConstraint} yields $\|\delta \locals\|\le \omega \Delta_i^L$. 
	Thus, due to the assumption $0<\omega\le\frac{N}{\|C_i\|}\Delta_i^L$,
	we obtain that 
	\begin{equation*}
		\|g_i-\tilde g_i\| = \|g_i-\left(g_i - C_i\cdot\sum_k  \Ik \delta \locals\right)\| \le N \|C_i\|\omega \Delta_i^L \le \Delta_i^L.
	\end{equation*}	
\end{proof}

\subsubsection{A Linear Recombination Approach}
\label{sec:linearApproach}

In this section, we consider a damping approach in order to satisfy Assumption~\eqref{eqn:assumptionDeltaS}. 
The damping parameter~$\alpha_i\in[0,1]$ is employed to linearly combine
the current gradient with the local corrections. This approach has the advantage that the local solution
process does not have to accept an initial Newton step as in \ref{eqn:localStart},
which might perhaps not exist or spoil the non-linear solution process.
We thus  do not state assumptions on the local solution process, 
but define the preconditioned gradient as
\begin{equation}
\label{eqn:linearUpdate}
	\tilde g_i = \alpha_i g_i - (1-\alpha_i)\cdot C_i \cdot \sum_k \Ik \locals,
\end{equation}
where $C_i$ is the inverse of the additive Schwarz preconditioner as defined in \eref{eqn:matrixC}. 
However, note that -- as we pointed out in Section~\ref{sec:preconditionedQuadraticModel} -- 
the computational cost for~$C_i$ and thus the preconditioned gradient~$\tilde g_i$ depends on the employed domain decomposition.

In order to compute a damping parameter $\alpha_i$ which satisfies \eref{eqn:assumptionDeltaS}, 
we estimate 
$$\|g_i-\tilde g_i\|=\left\|g_i-\left(\alpha_i g_i - (1-\alpha_i)\cdot C_i \cdot \sum_k \Ik \locals\right)\right\|$$ 
as follows
\begin{eqnarray*}
	\|g_i-\tilde g_i\| \le (1-\alpha_i) \left(\left\|g_i\right\|+\left\| C_i \cdot \sum_k \Ik \locals\right\|\right).
\end{eqnarray*}
Therefore, if $\alpha_i$ satisfies the inequality
\begin{eqnarray*}
	(1-\alpha_i) \left(\|g_i\|+ \|C_i \cdot \sum_k \Ik \locals\|\right) \le \Delta_i^L,	
\end{eqnarray*}
$\tilde g_i$ satisfies \eref{eqn:assumptionDeltaS}.
Thus, since $\alpha_i\in[0,1]$, we obtain that each $\tilde g_i$ given by \eref{eqn:linearUpdate} with
\begin{equation*}
	  \alpha_i = \min\left\{1, \max\left\{0, 1-\frac{\Delta_i^L}{\|g_i\|+ \|C_i \cdot \sum_k \Ik \locals\|}\right\}\right\}
\end{equation*}
satisfies \eref{eqn:assumptionDeltaS}.

Basically this update means that as long as $\Delta_i^L$ is sufficiently large we have $\tilde g_i=g_i^\text{ASPIN}$. 
But, if it turns out that a sufficient decrease cannot be achieved or if the
approximation strength of the preconditioned model is too small, and thus $\Delta_i^L$ is reduced, 
the original gradient is taken more and more into account. 
Furthermore, let us remark that the computation of $\alpha_i$ only depends on computable and known quantities.